\documentclass[sigconf]{acmart}
\usepackage[utf8]{inputenc} 
\usepackage[T1]{fontenc}    
\usepackage{hyperref}       
\usepackage{url}            
\usepackage{booktabs}       
\usepackage{amsfonts}       
\usepackage{nicefrac}       
\usepackage{microtype}      
\usepackage{epsfig}
\usepackage{graphicx}
\usepackage{amsmath}
\usepackage{amssymb}
\usepackage{amsthm}
\usepackage{array}

\usepackage{amsmath}
\usepackage{url}
\usepackage{algorithm}
\usepackage{diagbox}
\usepackage{bbm}
\usepackage{bm}

\usepackage{graphicx}
\usepackage{subcaption}
\usepackage{epstopdf}
\usepackage{color}
\usepackage{algorithm}
\usepackage{algorithmic}
\usepackage{multirow}
\usepackage{rotating}
\usepackage{alphalph}
\newcommand{\bbb}[1]{\boldsymbol{\mathbf{#1}}}
\def\beq{\begin{eqnarray}}
\def\eeq{\end{eqnarray}}
\def\noi{\noindent}
\def\nn{\nonumber}
\def\la{\langle}
\def\ra{\rangle}

\def\ghs{\hspace{0.003cm}} 
\newcommand{\fonttiny}{\fontsize{5}{5}\selectfont}

\def\fourfigwid{0.23\textwidth}

\def\objimgwid{0.9\textwidth}
\def\objimghei{0.085\textheight}
\def\E{\mathbb{E}}
\def\P{\mathcal{P}}
\usepackage{fancyhdr}
\usepackage{dsfont}
\usepackage{color}
\usepackage{setspace}

\newtheorem{lemma}{Lemma}
\newtheorem{theorem}{Theorem}
\newtheorem{definition}{Definition}
\newtheorem{proposition}{Proposition}

\newcommand{\fractt}{\tfrac}
\usepackage{multirow}
\usepackage{diagbox}

\def\cone{\textcolor[rgb]{0.9961,0,0}}
\def\ctwo{\textcolor[rgb]{0,0.5820,0.9}}
\def\cthree{\textcolor[rgb]{0,0.7,0}}
\definecolor{mygray}{gray}{.9}
\definecolor{mypink}{rgb}{.99,.91,.95}

\definecolor{mycyana}{RGB}{153,204,153}
\definecolor{mycyanb}{RGB}{255,153,102}
\definecolor{mycyanc}{RGB}{255,204,51}
\definecolor{mycyand}{RGB}{204,204,102}

\usepackage{textcomp,booktabs}
\usepackage{colortbl}

\usepackage[T1]{fontenc} \usepackage{ae} \usepackage{aecompl}

\copyrightyear{2020}
\acmYear{2020}
\setcopyright{acmcopyright}\acmConference[KDD '20]{Proceedings of the 26th ACM SIGKDD Conference on Knowledge Discovery and Data Mining}{August 23--27, 2020}{Virtual Event, CA, USA}
\acmBooktitle{Proceedings of the 26th ACM SIGKDD Conference on Knowledge Discovery and Data Mining (KDD '20), August 23--27, 2020, Virtual Event, CA, USA}
\acmPrice{15.00}
\acmDOI{10.1145/3394486.3403070}
\acmISBN{978-1-4503-7998-4/20/08}

\begin{document}

\fancyhead{}

\title{A Block Decomposition Algorithm for Sparse Optimization}

\author{Ganzhao Yuan$^{1}$, Li Shen$^{2}$, Wei-Shi Zheng$^{3,1}$}
\affiliation{$^1$Peng Cheng Laboratory, China \hspace{0.3cm}$^2$Tencent AI Lab, China\hspace{0.3cm}$^3$Sun Yat-sen University, China}

\affiliation{\small yuangzh@pcl.ac.cn,~ mathshenli@gmail.com, ~zhwshi@mail.sysu.edu.cn}

%
%
%
%

\begin{abstract}
Sparse optimization is a central problem in machine learning and computer vision. However, this problem is inherently NP-hard and thus difficult to solve in general. Combinatorial search methods find the global optimal solution but are confined to small-sized problems, while coordinate descent methods are efficient but often suffer from poor local minima. This paper considers a new block decomposition algorithm that combines the effectiveness of combinatorial search methods and the efficiency of coordinate descent methods. Specifically, we consider a random strategy or/and a greedy strategy to select a subset of coordinates as the working set, and then perform a global combinatorial search over the working set based on the original objective function. We show that our method finds stronger stationary points than Amir Beck et al.'s coordinate-wise optimization method. In addition, we establish the convergence rate of our algorithm. Our experiments on solving sparse regularized and sparsity constrained least squares optimization problems demonstrate that our method achieves state-of-the-art performance in terms of accuracy. For example, our method generally outperforms the well-known greedy pursuit method.


\end{abstract}



\maketitle

\section{Introduction}


This paper mainly focuses on the following nonconvex sparsity constrained / sparse regularized optimization problem:
\beq\label{eq:main}
\begin{split}
\boxed{ \min_{\bbb{x}}~f(\bbb{x}),~s.t.~\|\bbb{x}\|_0\leq s}\hspace{0.2cm}\text{or}\hspace{0.2cm}\boxed{\min_{\bbb{x}}~f(\bbb{x})+ \lambda \|\bbb{x}\|_0,}
\end{split}
\eeq
\noi where $\bbb{x} \in \mathbb{R}^n$, $\lambda$ is a positive constant, $s\in[n]$ is a positive integer, $f(\cdot)$ is assumed to be convex, and $\|\cdot\|_0$ is a function that counts the number of nonzero elements in a vector. Problem (\ref{eq:main}) can be rewritten as the following unified composite minimization problem (`$\triangleq$' means define):
\beq
&&~\min_{\bbb{x}}~F(\bbb{x})\triangleq  f(\bbb{x}) + h(\bbb{x}),~~~\text{with}~~~~~h(\bbb{x}) \triangleq h_{\text{cons}}~\text{or}~ h_{\text{regu}}.\text{~~~~~~~~~~~~} \nn
\eeq
\noi Here, $h_{\text{cons}}(\bbb{x}) \triangleq I_{\Psi}(\bbb{x}),~\Psi \triangleq \{\bbb{x}~|~\|\bbb{x}\|_0 \leq s\}$, $I_{\Psi}(\cdot)$ is an indicator function on the set $\Psi$ with $I_{\Psi}(\bbb{x})={\fonttiny \left\{
                                                                   \begin{array}{ll}
                                                                     0, &\hbox{$\bbb{x}\in\Psi$} \\
                                                                     \infty, &\hbox{$\bbb{x}\notin\Psi$}
                                                                   \end{array}
                                                                 \right.}$, and $h_{\text{regu}}(\bbb{x}) \triangleq \lambda \|\bbb{x}\|_0$. Problem (\ref{eq:main}) captures a variety of applications of interest in both machine learning and computer vision (e.g., sparse coding                                 \cite{Donoho06,aharon2006img,BaoJQS16}, sparse subspace clustering \cite{elhamifar2013sparse}).


This paper proposes a block decomposition algorithm using a proximal strategy and a combinatorial search strategy for solving the sparse optimization problem as in (\ref{eq:main}). We review existing methods in the literature and summarize the merits of our approach.

\textbf{$\blacktriangleright$ The Relaxed Approximation Method.} One popular method to solve Problem (\ref{eq:main}) is the convex or nonconvex relaxed approximation method \cite{candes2005decoding,zhang10a,xu2012l}. Many approaches such as $\ell_1$ norm, top-$k$ norm, Schatten $\ell_p$ norm, re-weighted $\ell_1$ norm, capped $\ell_1$ norm, and half quadratic function have been proposed for solving sparse optimization problems in the last decade. It is generally believed that nonconvex methods often achieve better accuracy than the convex counterparts \cite{YuanG19,Bi2014,yuan2016sparsempec}. However, minimizing the approximate function does not necessarily lead to the minimization of the original function in Problem (\ref{eq:main}). \textbf{$\spadesuit$} {Our method} directly controls the sparsity of the solution and minimize the original objective function.

%
%
%


\textbf{$\blacktriangleright$ The Greedy Pursuit Method.} This method is often used to solve sparsity constrained optimization problems. It greedily selects at each step one coordinate of the variables which have some desirable benefits \cite{tropp2007signal,dai2009subspace,needell2010signal,blumensath2008gradient,needell2009cosamp}. This method has a monotonically decreasing property and achieves optimality guarantees in some situations, but it is limited to solving problems with smooth objective functions (typically the square function). Furthermore, the solutions must be initialized to zero and may cause divergence when being incorporated to solve the bilinear matrix factorization problem \cite{BaoJQS16}. \textbf{$\spadesuit$} {Our method} is a greedy coordinate descent algorithm without forcing the initial solution to zero.


\textbf{$\blacktriangleright$ The Combinatorial Search Method.} This method is typically concerned with NP-hard problems \cite{conforti2014integer}. A naive strategy is an exhaustive search which systematically enumerates all possible candidates for the solution and picks the best candidate corresponding to the lowest objective value. The cutting plane method solves the convex linear programming relaxation and adds linear constraints to drive the solution towards binary variables, while the branch-and-cut method performs branches and applies cuts at the nodes of the tree having a lower bound that is worse than the current solution. Although in some cases these two methods converge without much effort, in the worse case they end up solving all $2^n$ convex subproblems. \textbf{$\spadesuit$} {Our method} leverages the effectiveness of combinatorial search methods.

\bbb{$\blacktriangleright$ The Proximal Gradient Method.} Based on the current gradient $\nabla f(\bbb{x}^k)$, the proximal gradient method \cite{beck2013sparsity,lu2014iterative,jain2014iterative,patrascu2015efficient,PatrascuN15,BLUMENSATH2009265} iteratively performs a gradient update followed by a proximal operation: $\bbb{x}^{k+1}= \text{prox}(\bbb{x}^k - \beta \nabla f(\bbb{x}^k);\beta,h)$. Here the proximal operator $\text{prox}(\bbb{a};\beta,h) = \arg \min_{\bbb{x}} ~\tfrac{1}{2}\|\bbb{x}-\bbb{a}\|_2^2 + \beta h(\bbb{x})$ can be evaluated analytically, and $\beta={1}/{L}$ is the step size with $L$ being the Lipschitz constant. This method is closely related to (block) coordinate descent \cite{nesterov2012efficiency,xu2013block,razaviyayn2013unified,hong2013iteration,lu2015complexity} in the literature. Due to its simplicity, many strategies (e.g., variance reduction \cite{johnson2013accelerating,xiao2014proximal,LiZALH16}, asynchronous parallelism \cite{liu2015asynchronous,recht2011hogwild}, and non-uniform sampling \cite{zhang2016accelerated}) have been proposed to accelerate proximal gradient method. However, existing works use a scalar step size and solve a first-order majorization/surrogate function via closed-form updates. Since Problem (\ref{eq:main}) is nonconvex, such a simple majorization function may not necessarily be a good approximation. \textbf{$\spadesuit$} {Our method} significantly outperforms proximal gradient method and inherits its computational advantages.






\textbf{Contributions.} The contributions of this paper are three-fold. \bbb{(i)} Algorithmically, we introduce a novel block decomposition method for sparse optimization (See Section \ref{sect:proposed}). \bbb{(ii)} Theoretically, we establish the optimality hierarchy of our algorithm and show that it always finds stronger stationary points than existing methods (See Section \ref{sect:hierarchy}). Furthermore, we prove the convergence rate of our algorithm (See Section \ref{sect:convergence}). Additional discussions for our method is provided in Section \ref{sect:dis}. \bbb{(iii)} Empirically, we have conducted experiments on some sparse optimization tasks to show the superiority of our method (See Section \ref{sect:exp}).

\textbf{Notations.} All vectors are column vectors and superscript $\mathsf{T}$ denotes transpose. For any vector $\bbb{x}\in \mathbb{R}^n$ and any $i\in \{1,2,...,n\}$, we denote by $\bbb{x}_i$ the $i$-th component of $\bbb{x}$. The Euclidean inner product between $\bbb{x}$ and $\bbb{y}$ is denoted by $\la \bbb{x},\bbb{y} \ra$ or $\bbb{x}^{\intercal}\bbb{y}$. $\bbb{1}$ is an all-one column vector, and $\bbb{e}_i$ is a unit vector with a $1$ in the $i$th entry and $0$ in all other entries. When $\beta$ is a constant, $\beta^t$ denotes the $t$-th power of $\beta$, and when $\beta$ is an optimization variable, $\beta^t$ denotes the value of $\beta$ in the $t$-th iteration. The number of possible combinations choosing $k$ items from $n$ without repetition is denoted by $C_n^k$. For any $B \in \mathbb{N}^k$ containing $k$ unique integers selected from $\{1, 2, ..., n\}$, we define $\bar{B} \triangleq   \{1,2,...,n\}\setminus B$ and denote $\bbb{x}_B$ as the sub-vector of $\bbb{x}$ indexed by $B$.

\section{The Proposed Block Decomposition Algorithm} \label{sect:proposed}

This section presents our block decomposition algorithm for solving (\ref{eq:main}). Our algorithm is an iterative procedure. In every iteration, the index set of variables is separated into two sets $B$ and $\bar{B}$, where $B$ is the working set. We fix the variables corresponding to $\bar{B}$, while minimizing a sub-problem on variables corresponding to $B$. The proposed method is summarized in Algorithm \ref{algo:main}.

\begin{algorithm}[!h]
\caption{ {\bf The Proposed Block Decomposition Algorithm} }
\begin{algorithmic}[1]
  \STATE Input: the size of the working set $k\in[n]$, the proximal point parameter $\theta>0$, and an initial feasible solution $\bbb{x}^0$. Set $t=0$.
\WHILE{not converge}
\STATE (S1) Employ some strategy to find a working set $B$ of size $k$. We define $\bar{B} \triangleq   \{1,2,...,n\}\setminus B$.
  \STATE (S2) Solve the following subproblem~\underline{globally}~using combinatorial search methods:
  \beq \label{eq:subprob}
   \bbb{x}^{t+1} \Leftarrow \arg \min_{\bbb{z}}~f(\bbb{z}) +h(\bbb{z}) + \tfrac{\theta}{2} \|\bbb{z}-\bbb{x}^{t}\|^2,s.t.~\bbb{z}_{\bar{B}} = \bbb{x}^t_{\bar{B}}
\eeq
\vspace{-10pt}
\STATE (S3) Increment $t$ by 1
\ENDWHILE
\end{algorithmic}
\label{algo:main}
\end{algorithm}

At first glance, Algorithm \ref{algo:main} might seem to be merely a block coordinate descent algorithm \cite{tseng2009coordinate} applied to (\ref{eq:main}). However, it has some interesting properties that are worth commenting on.

\noi \bbb{$\blacktriangleright$ Two New Strategies.} \bbb{(i)} Instead of using majorization techniques for optimizing over the block of the variables, we consider minimizing the original objective function. Although the subproblem is NP-hard and admits no closed-form solution, we can use an exhaustive search to solve it exactly. \bbb{(ii)} We consider a proximal point strategy for the subproblem in (\ref{eq:subprob}). This is to guarantee sufficient descent condition for the optimization problem and global convergence of Algorithm \ref{algo:main} (refer to Proposition \ref{theorem:convergence}).

\noi \bbb{$\blacktriangleright$ Solving the Subproblem Globally.} The subproblem in (\ref{eq:subprob}) essentially contains $k$ unknown decision variables and can be solved exactly within sub-exponential time $\mathcal{O}(2^k)$. Using the variational reformulation of $\ell_0$ pseudo-norm \footnote{For all $\bbb{x} \in \mathbb{R}^n$ with $\|\bbb{x}\|_{\infty} \leq \rho$, it always holds that $\|\bbb{x}\|_0 = \min_{\bbb{v}} ~\la \bbb{1},\bbb{v} \ra,~s.t.~\bbb{v} \in \{0,1\}^n,~|\bbb{x}|\leq \rho \bbb{v}.$}, Problem (\ref{eq:subprob}) can be reformulated as a mixed-integer optimization problem and solved by some global optimization solvers such as `CPLEX' or `Gurobi'. For simplicity, we consider a simple exhaustive search (a.k.a. generate and test method) to solve it. Specifically, for every coordinate of the $k$-dimensional subproblem, it has two states, i.e., zero/nonzero. We systematically enumerate the full binary tree to obtain all possible candidate solutions and then pick the best one that leads to the lowest objective value as the optimal solution. 




\noi \bbb{$\blacktriangleright$ Finding the Working Set.} We observe that it contains $C_n^k$ possible combinations of choice for the working set. One may use a cyclic strategy to alternatingly select all the choices of the working set. However, past results show that the coordinate gradient method results in faster convergence when the working set is chosen in an arbitrary order \cite{hsieh2008dual} or in a greedy manner \cite{tseng2009coordinate,hsieh2011fast}. This inspires us to use a random strategy or a greedy strategy for finding the working set. We remark that the combination of the two strategies is preferred in practice.

\hspace{12pt} \boxed{\text{Random strategy}.} We uniformly select one combination (which contains $k$ coordinates) from the whole working set of size $C_n^k$. One good benefit of this strategy is that our algorithm is ensured to find a block-$k$ stationary point (discussed later) in expectation.

\hspace{12pt} \boxed{\text{Greedy strategy}.} Generally speaking, we pick the top-$k$ coordinates that lead to the greatest descent when one variable is changed and the rest variables are fixed based on the current solution $\bbb{x}^t$. We denote $Z\triangleq \{i: \bbb{x}^t_i=0\}$ and $\bar{Z}\triangleq \{j: \bbb{x}^t_j \neq 0\}$. For $Z$, we solve a one-variable subproblem to compute the possible decrease for all $i \in Z$ of $\bbb{x}^t$ when changing from zero to nonzero:
\beq
\textstyle \forall i =1,...,|Z|,~ \bbb{c}_{i} = \min_{\alpha} F(\bbb{x}^t + \alpha \bbb{e}_i) - F(\bbb{x}^t). \nn
  \eeq
  \noi For $\bar{Z}$, we compute the decrease for each coordinate $j \in \bar{Z}$ of $\bbb{x}^t$ when changing from nonzero to exactly zero:
  \beq
\textstyle \forall j =1,...,|\bar{Z}|,~\bbb{d}_{j} = F(\bbb{x}^t + \alpha \bbb{e}_j) - F(\bbb{x}^t),~\alpha = \bbb{x}^t_j. \nn
  \eeq
  \noi We sort the vectors $\bbb{c}$ and $\bbb{d}$ in increasing order and then pick the top-$k$ coordinates as the working set.

\section{Optimality Analysis}\label{sect:hierarchy}

This section provides an optimality analysis of our method. We assume that $f(\bbb{x})$ is a smooth convex function with its gradient being $L$-Lipschitz continuous. In the sequel, we present some necessary optimal conditions for (\ref{eq:main}). Since the block-$k$ optimality condition is novel in this paper, it is necessary to clarify its relations with existing optimality conditions formally. We use $\breve{\bbb{x}},~\grave{\bbb{x}}$, and $\bar{\bbb{x}}$ to denote an arbitrary basic stationary point, an $L$-stationary point, and a block-$k$ stationary point, respectively.

\begin{definition} \label{def:basic}
(Basic Stationary Point) A solution $\breve{\bbb{x}}$ is called a basic stationary point if the following holds. $h \triangleq h_{\text{cons}}:\breve{\bbb{x}} = \arg\min_{\bbb{y}} f(\bbb{y}),~s.t.~|\bar{Z}|\leq k,~\bbb{y}_Z = \bbb{0}$; $h \triangleq h_{\text{regu}}: \breve{\bbb{x}} = \arg \min_{\bbb{y}}~f(\bbb{y}),~s.t.~\bbb{y}_Z = \bbb{0}$. Here, $Z\triangleq \{i| \breve{\bbb{x}}_i = 0\}$,~$\bar{Z}\triangleq \{j| \breve{\bbb{x}}_j \neq 0\}$.
\end{definition}

\noi \textbf{Remarks.} The basic stationary point states that the solution achieves its global optimality when the support set is restricted. The number of basic stationary points for $h \triangleq h_{\text{cons}}$ and $h \triangleq h_{\text{regu}}$ is $\sum_{i=0}^k C_n^i$ and $\sum_{i=0}^n C_n^i$, respectively. One good feature of the basic stationary condition is that the solution set is enumerable, which makes it possible to validate whether a solution is optimal for the original sparse optimization problem.

\begin{definition} \label{def:L:stationary}
($L$-Stationary Point) A solution $\grave{\bbb{x}}$ is an $L$-stationary point if it holds that: $\grave{\bbb{x}} = \arg\min_{\bbb{y}}~g(\bbb{y},\grave{\bbb{x}})+ h(\bbb{y})$ with $g(\bbb{y},\bbb{x})\triangleq f(\bbb{x}) + \la \nabla f(\bbb{x}),\bbb{y} - \bbb{x} \ra + \fractt{L}{2}\|\bbb{y} - \bbb{x}\|_2^2$.
\end{definition}

\noi \textbf{Remarks.} This is the well-known proximal thresholding operator. The term $g(\bbb{y},\bbb{x})$ is a majorization function of $f(\bbb{y})$ and it always holds that $f(\bbb{y})\leq  g(\bbb{y},\bbb{x})$ for all $\bbb{x}$ and $\bbb{y}$. Although it has a closed-form solution, this simple surrogate function may not be a good majorization/surrogate function for the non-convex problem.


\begin{definition} \label{def:block:k}
 (Block-$k$ Stationary Point) A solution $\bar{\bbb{x}}$ is a block-$k$ stationary point if it holds that:
\beq
\bar{\bbb{x}}\in \arg\min_{\bbb{z}}\P(\bbb{z};\bar{\bbb{x}},B)\triangleq \{F(\bbb{z}),s.t.\bbb{z}_{\bar{B}}=\bar{\bbb{x}}_{\bar{B}}\},\forall |B|=k.
\eeq
\end{definition}

\noi\textbf{Remarks.} \bbb{(i)} The concept of the block-$k$ stationary point is novel in this paper. Our method can inherently better explore the second-order / curvature information of the objective function. \bbb{(ii)} The sub-problem $\min_{\bbb{z}}~\P(\bbb{z};\bar{\bbb{x}},B)$ is NP-hard, and it takes sub-exponential time $\mathcal{O}(2^k)$ to solve it. However, since $k$ is often very small, it can be tackled by some practical global optimization methods. \bbb{(iii)} Testing whether a solution $\bar{\bbb{x}}$ is a block-$k$ stationary point deterministically requires solving $C_n^k$ subproblems, therefore leading to a total time complexity of $C_n^k \times \mathcal{O}(2^k)$. However, using a random strategy for finding the working set $B$ from $C_n^k$ combinations, we can test whether a solution $\bar{\bbb{x}}$ is the block-$k$ stationary point in expectation within a time complexity of $T\times \mathcal{O}(2^k)$ with the constant $T$ being the number of times which is related to the confidence of the probability.


\begin{table*}[!t]
\fontsize{5}{6}\selectfont
\centering
\scalebox{1.15}{\begin{tabular}{|p{1.2cm}|p{1.2cm}|p{1.2cm}|p{1.4cm}|p{1.4cm}|p{1.4cm}|p{1.4cm}|p{1.4cm}|p{1.4cm}|}
  \hline
 & {  \text{Basic} \text{Stat.}} & { L-\text{Stat.}} & {  \text{Block}-1 \text{Stat.}} & { \text{Block}-2 \text{Stat.} }&  { \text{Block}-3 \text{Stat.}}& {  \text{Block}-4 \text{Stat.}}&  { \text{Block}-5 \text{Stat.}}&  { \text{Block}-6 \text{Stat.}}\\
  \hline
$h =h_{\text{cons}}$ &  57 & 14 &  -- & 2 &  1 &  1&  1&  1 \\
  \hline
$h = h_{\text{regu}}$ &  64 & 56 &  9 & 3 &  1 &  1&  1&  1 \\
  \hline
\end{tabular}}
\caption{Number of points satisfying optimality conditions.}\label{tab:optimality}
\vspace{-10pt}
\end{table*}

The following proposition states the relations between the three types of the stationary point.
\begin{proposition} \label{proposition:hierarchy}
\textbf{Optimality Hierarchy between the Optimality Conditions.} The following relationship holds:\\ $\boxed{\text{Basic Stat. Point}} \overset{(a)}{\Leftarrow} \boxed{ \text{$L$-Stat. Point}} \overset{(b)}{\Leftarrow}  \boxed{  \text{Block-$k$ Stat. Point}}$\\$ \overset{(c)}{\Leftarrow}  \boxed{  \text{Block-$(k+1)$ Stat. Point}}   \overset{}{\Leftarrow}  ... \overset{}{\Leftarrow}  \boxed{  \text{Block-$n$ Stat. Point}}   \overset{(d)}{\Leftrightarrow}  \boxed{ \text{Optimal Point}}$ for sparse regularized (resp., for sparsity constrained) optimization problems with $k\geq 1$ (resp., $k\geq 2$).



\begin{proof}
We denote $\Gamma_s(\bbb{x})$ as the operator that sets all but the largest (in magnitude) $s$ elements of $\bbb{x}$ to zero.

\bbb{(a)} First, we prove that an $L$-stationary point $\grave{\bbb{x}}$ is also a basic stationary point $\breve{\bbb{x}}$ when $h \triangleq h_{\text{cons}}$. For an $L$-stationary point, we have $\grave{\bbb{x}} = \Gamma_s(\grave{\bbb{x}} - (\nabla f(\grave{\bbb{x}}))/L)$. This implies that there exists an index set $S$ such that $\breve{\bbb{x}}_S = \breve{\bbb{x}}_S - (\nabla f(\breve{\bbb{x}}))_S/L$ and $\breve{\bbb{x}}_{\{1,...,n\}\setminus S} = \bbb{0}$, which is the optimal condition for a basic stationary point.

Second, we prove that an $L$-stationary point $\grave{\bbb{x}}$ is also a basic stationary point $\breve{\bbb{x}}$ when $h \triangleq h_{\text{regu}}$. Using Definition \ref{def:L:stationary}, we have the following closed-form solution for $\grave{\bbb{x}}$:
{\small \beq
\grave{\bbb{x}}_i = \left\{
              \begin{array}{ll}
                \grave{\bbb{x}}_i - \nabla_i f(\grave{\bbb{x}})/L, & \hbox{$(\grave{\bbb{x}}_i - \nabla_i f(\grave{\bbb{x}})/L)^2 > 2 \lambda/L$;} \\
                0, & \hbox{$\text{else}.$}
              \end{array}
            \right.\nn
 \eeq }{\normalsize}\noi This implies that there exists a support set $S$ such that $\breve{\bbb{x}}_S = \breve{\bbb{x}}_S - (\nabla f(\breve{\bbb{x}}))_S/L$, which is the optimal condition for a basic stationary point. Defining $Z\triangleq \{i| \grave{\bbb{x}}_i = 0\}$,~$\bar{Z}\triangleq \{j| \grave{\bbb{x}}_j \neq 0\}$, we note that $\forall i\in Z,~|\nabla f(\grave{\bbb{x}})|_i \leq \sqrt{2 \lambda L}$, and $\forall j \in \bar{Z},~(\nabla f(\grave{\bbb{x}}))_j = 0,~|\grave{\bbb{x}}_j | \geq \min(\sqrt{2 \lambda/L})$.

\bbb{(b)} First, we prove that a block-2 stationary point is also an $L$-stationary point for $h \triangleq h_{\text{cons}}$. Given a vector $\bbb{a}\in \mathbb{R}^n$, we consider the following optimization problem:
\beq \label{eq:CNK}
\bbb{z}^*_B= \arg \min_{\bbb{z}_B,~\forall |B|=2}~\|\bbb{z}-\bbb{a}\|_2^2,~s.t.~\|\bbb{z}_B\|_0 + \|\bbb{z}_{{\bar{B}}}\|_0 \leq s,
\eeq
\noi which essentially contains $C_n^2$ 2-dimensional subproblems. It is not hard to validate that (\ref{eq:CNK}) achieves the optimal solution with $\bbb{z}^* = \Gamma_s(\bbb{a})$. For any block-2 stationary point $\bar{\bbb{x}}$, we have $\bar{\bbb{x}}_B= \arg \min_{\bbb{z}_B}~\|\bbb{z}-(\bar{\bbb{x}} - \nabla f(\bar{\bbb{x}})/L)\|_2^2,~s.t.~\|\bbb{z}_B\|_0 + \|\bbb{z}_{\bar{B}}\|_0 \leq s$. Applying this conclusion with $\bbb{a}=\bar{\bbb{x}}-\nabla f(\bar{\bbb{x}})/L$, we have $\bar{\bbb{x}}=\Gamma_s(\bar{\bbb{x}}-\nabla f(\bar{\bbb{x}})/L)$.

Second, we prove that a block-1 stationary point is also an $L$-stationary point for $h \triangleq h_{\text{regu}}$. Assume that the convex objective function $f(\cdot)$ has coordinate-wise Lipschitz continuous gradient with constant $\bbb{s}_i,~\forall i=1,2,...,n$. For all $\bbb{x}\in\mathbb{R}^n,~t\in\mathbb{R},~i=1,2,...n$, it holds that \cite{nesterov2012efficiency}: $ f(\bbb{x}+ t \bbb{e}_i) \leq Q_i(\bbb{x},t)\triangleq f(\bbb{x}) + \la \nabla f(\bbb{x}),t\bbb{e}_i \ra + \fractt{\bbb{s}_i}{2}\| t\bbb{e}_i\|_2^2.$ Any block-1 stationary point must satisfy the following relation: $0 \in \arg \min_{t}~Q_i(\bar{\bbb{x}},t)+ \lambda \|\bar{\bbb{x}}_i + t \|_0,~\forall i.$ We have the following optimal condition for $\bar{\bbb{x}}$ with $k=1$: ${\tiny \bar{\bbb{x}}_i = {  \left\{
              \begin{array}{ll}
                (\bar{\bbb{x}}_i - \nabla_i f(\bar{\bbb{x}})/\bbb{s}_i  ), & \hbox{$(\bar{\bbb{x}}_i - \nabla_i f(\bar{\bbb{x}})/\bbb{s}_{i})^2 > {2\lambda}/{\bbb{s}_i}$;} \\
                0, & \hbox{$\text{else}$.}
              \end{array}
            \right\}}.}$ Since $\forall i,~\bbb{s}_i \leq L$, the latter formulation implies the former one.

%
%
%
%


\bbb{(c)} Assume $k_1\geq k_2$. The subproblem for the block-$k_2$ stationary point is a subset of those of the block-$k_1$ stationary point. Therefore, the block-$k_1$ stationary point implies the block-$k_2$ stationary point.

\bbb{(d)} Obvious.

\end{proof}

\end{proposition}

\noi\textbf{Remarks.} It is worthwhile to point out that the seminal works of Amir Beck
et al. present a coordinate-wise optimality condition for sparse optimization \cite{beck2013sparsity,Beck2016,BeckH19,BeckH20}. However, our block-$k$ condition is stronger since their optimality condition corresponds to $k=1$ in our optimality condition framework.

\vspace{0.5cm}

\textbf{A Running Example.} We consider the following sparsity constrained / sparse regularized optimization problem: $\min_{\bbb{x}\in\mathbb{R}^n}~\tfrac{1}{2}\bbb{x}^{\intercal}\bbb{Qx} + \bbb{x}^{\intercal}\bbb{p} + h(\bbb{x})$. Here, $n=6$,~$\bbb{Q}=\bbb{c}\bbb{c}^{\intercal}+\bbb{I}$,~$\bbb{p}=\bbb{1}$,~$\bbb{c}=[1~2~3~4~5~6]^{\intercal}$. The parameters for $h_{\text{cons}}(\bbb{x})$ and $h_{\text{regu}}(\bbb{x})$ are set to $(s,~\lambda)=(4,~0.01)$. The stationary point distribution of this example can be found in Table \ref{tab:optimality}. This problem contains $\sum_{i=0}^4 C_6^i = 57$ basic stationary points for $h \triangleq h_{\text{cons}}$, while it has $\sum_{i=0}^6 C_6^i=2^6 = 64$ basic stationary points for $h \triangleq h_{\text{regu}}$. As $k$ becomes large, the newly introduced type of local minimizer (i.e., block-$k$ stationary point) becomes more restricted in the sense that it has a smaller number of stationary points. Moreover, any block-3 stationary point is also the unique global optimal solution for this example.

\section{Convergence Analysis} \label{sect:convergence}

This section provides some convergence analysis for Algorithm \ref{algo:main}. We assume that $f(\bbb{x})$ is a smooth convex function with its gradient being $L$-Lipschitz continuous, and the working set of size $k$ is selected randomly and uniformly (sample with replacement). Due to space limitations, some proofs are placed into the \bbb{Appendix}.

%





\begin{proposition} \label{theorem:convergence}

\textbf{Global Convergence.} Letting $\{\bbb{x}^t\}_{t=0}^{\infty}$ be the sequence generated by Algorithm \ref{algo:main}, we have the following results. \bbb{(a)} It holds that: $F(\bbb{x}^{t+1})   \leq F(\bbb{x}^t) - \tfrac{\theta}{2}\|\bbb{x}^{t+1}-\bbb{x}^{t}\|^2,~\lim_{t\rightarrow \infty} \E[\|\bbb{x}^{t+1} - \bbb{x}^t\|] = 0$. \bbb{(b)} As $t\rightarrow \infty$, $\bbb{x}^t$ converges to the block-$k$ stationary point $\bar{\bbb{x}}$ of (\ref{eq:main}) in expectation.


\end{proposition}

\noi \textbf{Remarks.} Coordinate descent may cycle indefinitely if each minimization step contains multiple solutions \cite{Powell1973}. The introduction of the proximal point parameter $\theta>0$ is necessary for our nonconvex problem since it guarantees sufficient decrease condition, which is essential for global convergence. Our algorithm is guaranteed to find the block-$k$ stationary point, but it is in expectation.

We prove the convergence rate of our algorithm for sparsity constrained optimization with $h \triangleq h_{\text{cons}}$.

\begin{theorem}


\textbf{Convergence Rate for Sparsity Constrained Optimization.} Assume that $f(\cdot)$ is $\sigma$-strongly convex, and Lipschitz continuous such that $\forall t,~\|\nabla f(\bbb{x}^t)\|_{2}^2\leq \tau$ for some positive constant $\tau$. Denoting $\alpha \triangleq {\fractt{n\theta}{k \sigma}}/(1+\fractt{n\theta}{k\sigma})$, we have the following results:
$$\E[F(\bbb{x}^{t}) - F(\bar{\bbb{x}})] \leq   (F(\bbb{x}^{0}) - F(\bar{\bbb{x}}))\alpha^t +  \tfrac{\tau}{2\theta} \tfrac{\alpha}{1-\alpha},$$
$$\E[\tfrac{\sigma}{4}\|\bbb{x}^{t+1} - \bar{\bbb{x}}\|^2_2] \leq   \tfrac{2 n   \theta}{k}      (F(\bbb{x}^{0}) - F(\bar{\bbb{x}}))\alpha^t + \tfrac{n}{k}\fractt{\tau}{1-\alpha}.$$

\begin{proof}
\bbb{(a)} First of all, we define the zero set and nonzero set of the solution $\bbb{x}^{t+1}$ as follows:
\beq
S\triangleq \{i~|~i\in B,~\bbb{x}_{i}^{t+1} \neq 0\},~Q\triangleq \{i~|~i\in B,~\bbb{x}_{i}^{t+1} = 0\}.\nn
\eeq
\noi Using the optimality of $\bbb{x}^{t+1}$ for the subproblem, we obtain
\beq \label{eq:optimality:xk1}
(\nabla f(\bbb{x}^{t+1}))_{S} + \theta (\bbb{x}^{t+1}_S-\bbb{x}^{t}_S) = 0 
\eeq

\noi We derive the following inequalities:
\beq \label{eq:strongly:convex:sparse:rate}
&&  \E[f(\bbb{x}^{t+1}) - f(\bar{\bbb{x}})] \nn\\
&\overset{(a)}{\leq}&\E[ \la \bbb{x}^{t+1} - \bar{\bbb{x}},~\nabla f(\bbb{x}^{t+1}) \ra - \fractt{\sigma}{2}\|\bbb{x}^{t+1} - \bar{\bbb{x}}\|_2^2 ]\nn\\
&\overset{(b)}{=} & \fractt{n}{k}  \E[ \la \bbb{x}_{B}^{t+1} - \bar{\bbb{x}}_{B},~(\nabla f(\bbb{x}^{t+1}))_{B}\ra - \fractt{\sigma}{2}\|\bbb{x}^{t+1}_B - \bar{\bbb{x}}_B\|_{2}^2 ]   \nn\\
&\overset{(c)}{\leq} & \fractt{n}{k}  \fractt{\sigma}{2}  \E[\|(\nabla f(\bbb{x}^{t+1}))_{B}/\sigma\|_{2}^2 ]  \nn\\
&\overset{(d)}{=} & \fractt{n}{k}  \fractt{1}{2\sigma}  \left( \E[\|(\nabla f(\bbb{x}^{t+1}))_{S}\|_{2}^2 ] + \E[\|(\nabla f(\bbb{x}^{t+1}))_{Q}\|_{2}^2 ]\right)  \nn\\
& \overset{(e)}{\leq} & \E[ \fractt{n}{k}\fractt{1}{2\sigma}  [ \| \theta(\bbb{x}_S^t - \bbb{x}^{t+1}_S)\|_2^2  +\|(\nabla f(\bbb{x}^{t+1}))_Q\|_2^2 ] ]\nn\\
& \overset{(f)}{\leq} &   \E[ \fractt{n}{k}\fractt{\theta^2}{2\sigma}\| \bbb{x}^t - \bbb{x}^{t+1}\|_2^2  +  \fractt{n \tau}{2\sigma k}  ]\nn\\
& \overset{(g)}{\leq} & \E[   \fractt{n\theta}{\sigma k} [f(\bbb{x}^t)-f(\bbb{x}^{t+1})]  + \fractt{n \tau}{2\sigma k} ]\nn\\
& \overset{(h)}{=} & \E[ \fractt{n\theta }{\sigma k}[(f(\bbb{x}^t)-f(\bar{\bbb{x}}))-(f(\bbb{x}^{t+1})-f(\bar{\bbb{x}}))]  + \fractt{n \tau}{ 2\sigma k} ]~~~~~
\eeq
\noi where step $(a)$ uses the strongly convexity of $f(\cdot)$; step $(b)$ uses the fact that the working set $B$ is selected with $\fractt{k}{n}$ probability; step $(c)$ uses the inequality that $\la \bbb{x},~\bbb{a}\ra - \fractt{ \sigma}{2}\|\bbb{x}\|_2^2 = \fractt{\sigma}{2}\|\bbb{a}/\sigma\|_2^2-\fractt{\sigma}{2}\|\bbb{x}-\bbb{a}/\sigma\|_2^2\leq\fractt{\sigma}{2}\|\bbb{a}/\sigma\|_2^2$ for all $\bbb{a},~\bbb{x}$; step $(d)$ uses the fact that $B = S\cup Q$, step $(e)$ uses (\ref{eq:optimality:xk1}); step $(f)$ uses the fact that $\forall \bbb{x},~\|\bbb{x}_S\|_2^2\leq \|\bbb{x}\|_2^2$ and the Lipschitz continuity of $f(\cdot)$ that $\forall t,~\|\nabla f(\bbb{x}^{t+1})\|_2^2\leq \tau$; step $(g)$ uses the sufficient decrease condition that $\fractt{\theta}{2}\|\bbb{x}^{t+1}-\bbb{x}^{t}\|^2 \leq F(\bbb{x}^t) - F(\bbb{x}^{t+1})$; step $(h)$ uses $f(\bbb{x}^t)-f(\bbb{x}^{t+1})=(f(\bbb{x}^t)-f(\bar{\bbb{x}}))-(f(\bbb{x}^{t+1})-f(\bar{\bbb{x}}))$.

\noi From (\ref{eq:strongly:convex:sparse:rate}), we have the following inequalities:
\beq
&\E[(1+\fractt{n\theta}{ks}) (f(\bbb{x}^{t+1}) - f(\bar{\bbb{x}}))] \leq  \E[\fractt{n \theta}{ks}\cdot(f(\bbb{x}^t)-f(\bar{\bbb{x}})) + \fractt{n \tau}{2ks} ] \nn\\
&\E[f(\bbb{x}^{t+1}) - f(\bar{\bbb{x}})] \leq \E[ \alpha (f(\bbb{x}^t)-f(\bar{\bbb{x}})) + \fractt{\fractt{n}{2ks}}{\fractt{n \theta}{k \sigma}} \alpha \tau ]\nn\\
& \E[ f(\bbb{x}^{t+1}) - f(\bar{\bbb{x}})] \leq \E[ \alpha (f(\bbb{x}^t)-f(\bar{\bbb{x}})) +  \fractt{\alpha \tau}{2\theta}] \nn
\eeq

\noi Solving this recursive formulation, we have:
\beq \label{eq:convergence:rate:sparse:conc:1}
\E[ f(\bbb{x}^{t}) - f(\bar{\bbb{x}}) ] &\leq&  \E[ \alpha^t (f(\bbb{x}^{0}) - f(\bar{\bbb{x}}))] +  \tau \sum_{i=1}^t \alpha^i \nn\\
&=& \E[\alpha^t (f(\bbb{x}^{0}) - f(\bar{\bbb{x}}))] +  \fractt{\tau }{2\theta} \cdot \fractt{\alpha(1-\alpha^t)}{1-\alpha}\nn\\
&\leq & \E[\alpha^t (f(\bbb{x}^{0}) - f(\bar{\bbb{x}}))] +  \fractt{\tau }{2\theta} \cdot \fractt{\alpha}{1-\alpha}\nn
\eeq
\noi Since $\bbb{x}^t$ is always a feasible solution for all $t=1,2,...\infty$, we have $F(\bbb{x}^t)=f(\bbb{x}^t)$.

\noi \bbb{(b)} We now prove the second part of this theorem. First, we derive the following inequalities:
\beq \label{eq:bound:xt1xt}
\E[\|\bbb{x}^{t+1}-\bbb{x}^{t}\|_2^2] &\overset{(a)}{\leq} & \E[\fractt{2}{\theta} \left(F(\bbb{x}^t) - F(\bbb{x}^{t+1}) \right) ]\nn\\
&\overset{(b)}{\leq}& \E[\fractt{2}{\theta} \left( F(\bbb{x}^t) - F(\bar{\bbb{x}}) \right)] \nn\\
&\overset{(c)}{\leq}& \E[ \fractt{2}{\theta} \alpha^t (f(\bbb{x}^{0}) - f(\bar{\bbb{x}}))] +  \fractt{\tau \alpha }{\theta^2 (1-\alpha)}~~~~~~~~
\eeq
\noi where step $(a)$ uses the sufficient decrease condition; step $(b)$ uses the fact that $F(\bar{\bbb{x}}) \leq F({\bbb{x}}^{t+1})$; step $(c)$ uses the result in (\ref{eq:convergence:rate:sparse:conc:1}).

Second, we have the following results:
\beq \label{eq:strongly:convex:sparse:rate:eq}
 &&\E[\fractt{\sigma}{2}\|\bbb{x}^{t+1} - \bar{\bbb{x}}\|_2^2] \nn\\
 &\overset{(a)}{\leq}& \E[\la \bbb{x}^{t+1} - \bar{\bbb{x}},~\nabla f(\bbb{x}^{t+1}) \ra  + f(\bar{\bbb{x}})  - f(\bbb{x}^{t+1})]  \nn\\
&\overset{(b)}{\leq}& \E[ \la \bbb{x}^{t+1} - \bar{\bbb{x}},~\nabla f(\bbb{x}^{t+1}) \ra ] \nn\\
&\overset{(c)}{\leq}&  \E[\|\bbb{x}^{t+1} - \bar{\bbb{x}}\| \cdot \|\nabla f(\bbb{x}^{t+1})\|]
\eeq
\noi where step $(a)$ uses the strongly convexity of $f(\cdot)$; step $(b)$ uses the fact that $f(\bar{\bbb{x}}) \leq f({\bbb{x}}^{t+1})$; step $(c)$ uses the Cauchy-Schwarz inequality.

From (\ref{eq:strongly:convex:sparse:rate:eq}), we further have the following results:
\beq \label{eq:strongly:convex:sparse:eq2}
\E[\fractt{ \sigma}{4}\|\bbb{x}^{t+1} - \bar{\bbb{x}}\|^2_2] &\overset{(a)}{\leq}&\E[ \|\nabla f(\bbb{x}^{t+1})\|_2^2  ]\nn\\
&\overset{}{ = } & \fractt{n}{k} \E[  \|\nabla_B f(\bbb{x}^{t+1})\|_2^2  ]\nn\\
&\overset{(b)}{ = } &  \fractt{n}{k}  \E[  \|\nabla_S f(\bbb{x}^{t+1})\|_2^2 +  \|\nabla_Q f(\bbb{x}^{t+1})\|_2^2   ]\nn\\
&\overset{(c)}{ \leq } &   \fractt{n}{k} \E[  \theta^2\| \bbb{x}^{t+1}_S  - \bbb{x}^{t}_S \|_2^2] +   \fractt{n}{k} \tau   \nn\\
&\overset{(d)}{ = } &   \fractt{n}{k}  \E[ {2 \theta} \alpha^t (f(\bbb{x}^{0}) - f(\bar{\bbb{x}})) +   \fractt{\tau}{1-\alpha} ]  \nn
\eeq
\noi where step $(a)$ uses the strongly convexity of $f(\cdot)$; step $(b)$ uses the fact that $B = S \cup Q$; step $(c)$ uses the assumption that $\|\nabla f(\bbb{x}^t)\|_2^2\leq \tau$ for all $\bbb{x}^t$ and the optimality of $\bbb{x}^{t+1}$ in (\ref{eq:optimality:xk1}); step $(d)$ uses (\ref{eq:bound:xt1xt}). Therefore, we finish the proof of this theorem.

\end{proof}
\end{theorem}

\noi \textbf{Remarks.} Our results of convergence rate are similar to those of the gradient hard thresholding pursuit as in \cite{YuanLZ17}. The first term and the second term for our convergence rate are called parameter estimation error and statistical error, respectively. While their analysis relies on the conditions of restricted strong convexity/smoothness, our study relies on the requirements of generally strong convexity/smoothness.


We prove the convergence rate of our algorithm for sparse regularized optimization with $h\triangleq h_{\text{regu}}$.

\begin{theorem} \label{theorem:sparse:regularized:1}
\textbf{Convergence Rate for Sparse Regularized Optimization.} For $h\triangleq h_{\text{regu}}$, we have the following results:

\bbb{(a)} It holds that $\forall i,~|\bbb{x}_i^{t}|\geq \delta>0$ with $\bbb{\bbb{x}}_i^t\neq 0$. Whenever $\bbb{x}^{t+1} \neq \bbb{x}^t$, we have $\|\bbb{x}^{t+1}-\bbb{x}^t\|_2^2 \geq \fractt{k\delta^2}{n}$ and the objective value is decreased at least by $D$. The solution changes at most ${\bar{J}}$ times in expectation for finding a block-$k$ stationary point $\bar{\bbb{x}}$. Here $\delta$, $D$, and ${\bar{J}}$ are respectively defined as:
\beq \label{eq:DD}
\delta\triangleq \min(\sqrt{ \tfrac{2\lambda}{\theta + L} },\min(|\bbb{x}^0|)),~D \triangleq \tfrac{k \theta\delta^2}{2n},~{\bar{J}} \triangleq \tfrac{F(\bbb{x}^0) - F(\bar{\bbb{x}})}{D}.
\eeq

\bbb{(b)} Assume that $f(\cdot)$ is generally convex, and the solution is always bounded with $\|\bbb{x}^t\|_{\infty}\leq \rho,~\forall t$. If the support set of $\bbb{x}^t$ does not changes for all $t=0,1,...,\infty$, Algorithm \ref{algo:main} takes at most $V_1$ iterations in expectation to converge to a stationary point $\bar{\bbb{x}}$ satisfying $F(\bbb{x}^t)-F(\bar{\bbb{x}}) \leq \epsilon$. Moreover, Algorithm \ref{algo:main} takes at most $V_1 \times {\bar{J}}$ iterations in expectation to converge to a stationary point $\bar{\bbb{x}}$ satisfying $F(\bbb{x}^t)-F(\bar{\bbb{x}}) \leq \epsilon$. Here, $V_1$ is defined as:
\beq \label{eq:convex:boundt}
V_1 = {\max( \tfrac{4 \nu^2}{\theta} , \sqrt{  \tfrac{2 \nu^2 (F(\bbb{x}^0) - F(\bar{\bbb{x}}))}{\theta } })}/{\epsilon },~\text{with}~\nu\triangleq \tfrac{2n \rho \sqrt{k}\theta}{k}.
\eeq

\bbb{(c)} Assume that $f(\cdot)$ is $\sigma$-strongly convex. If the support set of $\bbb{x}^t$ does not changes for all $t=0,1,...,\infty$, Algorithm \ref{algo:main} takes at most $V_2$ iterations in expectation to converge to a stationary point $\bar{\bbb{x}}$ satisfying $F(\bbb{x}^t)-F(\bar{\bbb{x}}) \leq \epsilon$. Moreover, Algorithm \ref{algo:main} takes at most $V_2 \times {\bar{J}}$ iterations in expectation to converge to a stationary point $\bar{\bbb{x}}$ satisfying $F(\bbb{x}^t)-F(\bar{\bbb{x}}) \leq \epsilon$. Here, $V_2$ is defined as:
\beq \label{eq:convex:boundt2}
V_2 = \log_{\alpha} ({\epsilon}/{(F(\bbb{x}^{0}) - F(\bar{\bbb{x}}))}),~\text{with}~\alpha\triangleq {\tfrac{n\theta}{k \sigma}}/(1+\tfrac{n\theta}{k\sigma}).
\eeq

\noi

\end{theorem}


\noi \textbf{Remarks.} \bbb{(i)} When the support set is fixed, the optimization problem reduces to a convex problem. \bbb{(ii)} We derive a upper bound for the number of changes $\bar{J}$ for the support set in \bbb{(a)}, and a upper bound on the number of iterations $V_1$ (or $V_2$) performed after the support set is fixed in \bbb{(b)} (or \bbb{(c)}). Multiplying these two bounds, we can establish the upper bound of the number of iterations for Algorithm \ref{algo:main} to converge. However, these bounds are not tight enough. 

The following theorem establishes an improved convergence rate of our algorithm with $h \triangleq h_{\text{regu}}$.

\begin{theorem}\label{theorem:sparse:regularized:2}
\textbf{Improved Convergence Rate for Sparse Regularized Optimization.} For $h\triangleq h_{\text{regu}}$, we have the following results:

\bbb{(a)} Assume that $f(\cdot)$ is generally convex, and the solution is always bounded with $\|\bbb{x}^t\|_{\infty}\leq \rho,~\forall t$. Algorithm \ref{algo:main} takes at most $N_1$ iterations in expectation to converge to a block-$k$ stationary point $\bar{\bbb{x}}$ satisfying $F(\bbb{x}^t)-F(\bar{\bbb{x}})\leq \epsilon$, where $N_1= ( \fractt{\bar{J}}{D} + \fractt{1}{\epsilon})  \times \max( \fractt{4 \nu^2}{\theta} , \sqrt{  \fractt{2 \nu^2 (F(\bbb{x}^0) - F(\bar{\bbb{x}}) - D)}{\theta } })   $.

\bbb{(b)} Assume that $f(\cdot)$ is $\sigma$-strongly convex. Algorithm \ref{algo:main} takes at most $N_2$ iterations in expectation to converge to a block-$k$ stationary point $\bar{\bbb{x}}$ satisfying $F(\bbb{x}^t)-F(\bar{\bbb{x}})\leq \epsilon$, where $N_2=  \bar{J} \log_{\alpha} ( \fractt{D}{(F(\bbb{x}^0)-F(\bar{\bbb{x}}))} )  + \log_{\alpha} (\fractt{\epsilon}{F(\bbb{x}^{0}) - D - F(\bar{\bbb{x}})})$.

\end{theorem}

\noi \textbf{Remarks.} \bbb{(i)} Our proof of Theorem \ref{theorem:sparse:regularized:2} is based on the results in Theorem \ref{theorem:sparse:regularized:1} and a similar iterative bounding technique as in \cite{lu2014iterative}. \bbb{(ii)} If $\bar{J}\geq 2$ and the accuracy $\epsilon$ is sufficiently small such that $\epsilon\leq \fractt{D}{2}$, we have  $\fractt{\bar{J}}{D} + \fractt{1}{\epsilon} \leq \fractt{\bar{J}}{2 \epsilon} + \fractt{1}{\epsilon} \leq  \fractt{\bar{J}}{2 \epsilon} + \fractt{ \bar{J} / 2}{\epsilon} = \tfrac{\bar{J}}{\epsilon}$, leading to $ (\tfrac{\bar{J}}{D} + \fractt{1}{\epsilon})\times  \max( \fractt{4 \nu^2}{\theta} , \sqrt{ {2 \nu^2 (F(\bbb{x}^0) - F(\bar{\bbb{x}}) - D)}/{\theta } })   \leq \frac{\bar{J}}{\epsilon}\times \max( \tfrac{4 \nu^2}{\theta} , \sqrt{  {2 \nu^2 (F(\bbb{x}^0) - F(\bar{\bbb{x}}) )}/{\theta } })$. Using the same assumption and strategy, we have $\bar{J} \log_{\alpha} ({D}/{(F(\bbb{x}^0)-F(\bar{\bbb{x}}))} )  + \log_{\alpha} ({\epsilon}/{(F(\bbb{x}^{0}) - D - F(\bar{\bbb{x}}))}) \leq  \bar{J}\times \log_{\alpha} ({\epsilon}/{(F(\bbb{x}^{0}) - F(\bar{\bbb{x}}))})$. In this situation, the bounds in Theorem \ref{theorem:sparse:regularized:2} are tighter than those in Theorem \ref{theorem:sparse:regularized:1}.

\begin{figure*} [!t]
\centering
      \begin{subfigure}{\fourfigwid}\includegraphics[width=\objimgwid,height=\objimghei]{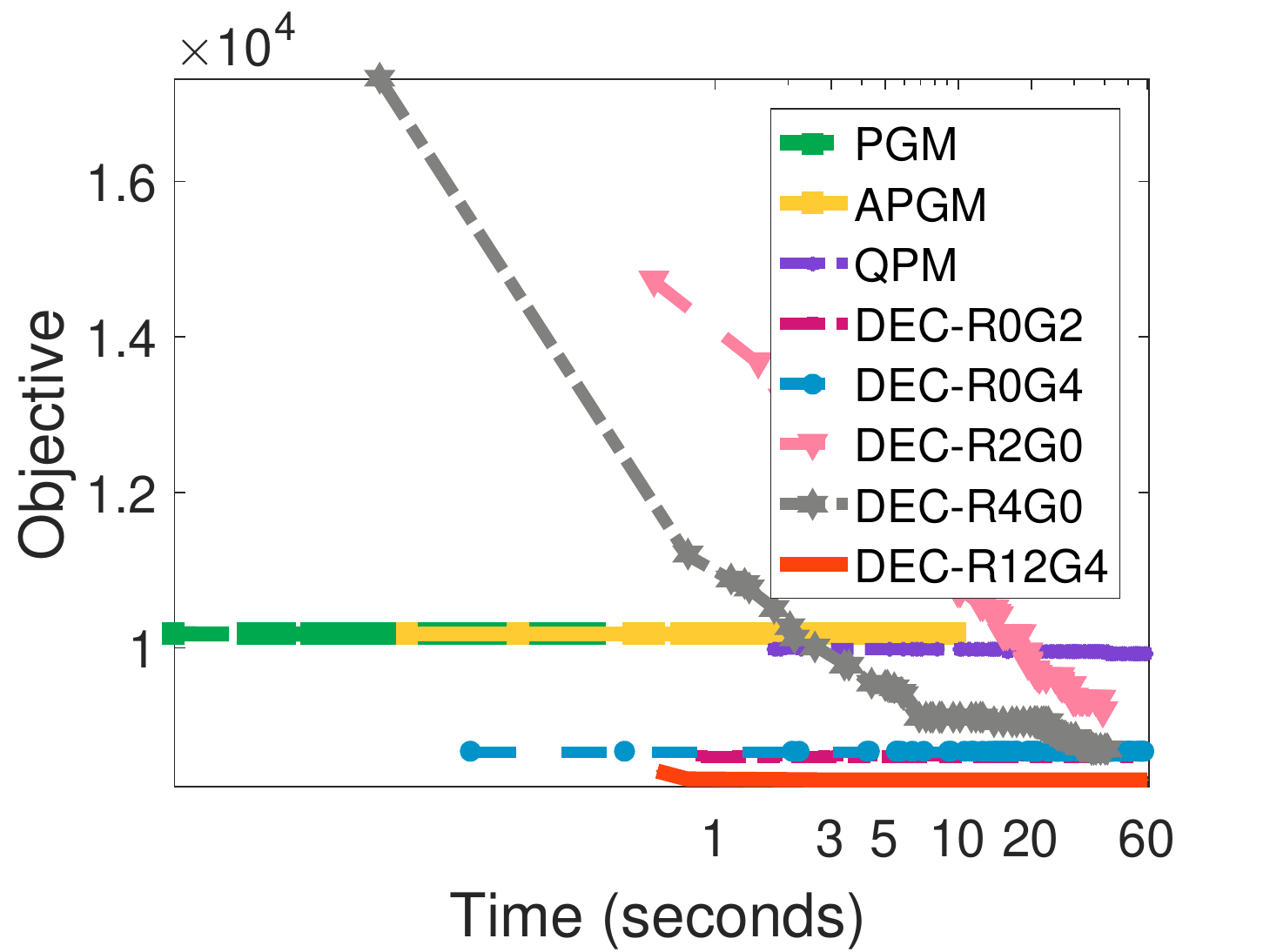}\vspace{-6pt} \caption{\scriptsize s=20 on random-256-1024 }\end{subfigure}\ghs
     \begin{subfigure}{\fourfigwid}\includegraphics[width=\objimgwid,height=\objimghei]{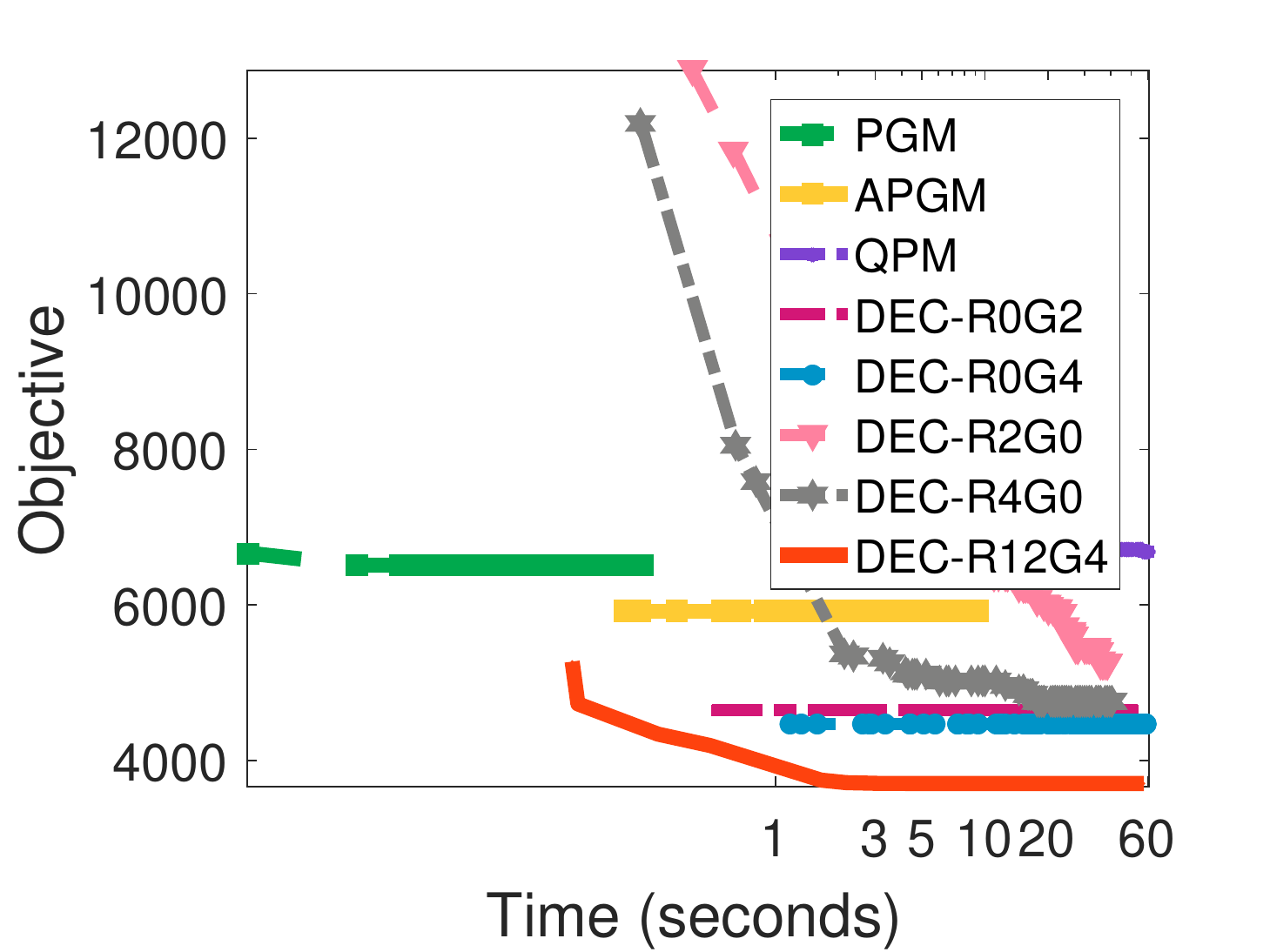}\vspace{-6pt} \caption{ \scriptsize s=40 on random-256-1024} \end{subfigure}\ghs
      \begin{subfigure}{\fourfigwid}\includegraphics[width=\objimgwid,height=\objimghei]{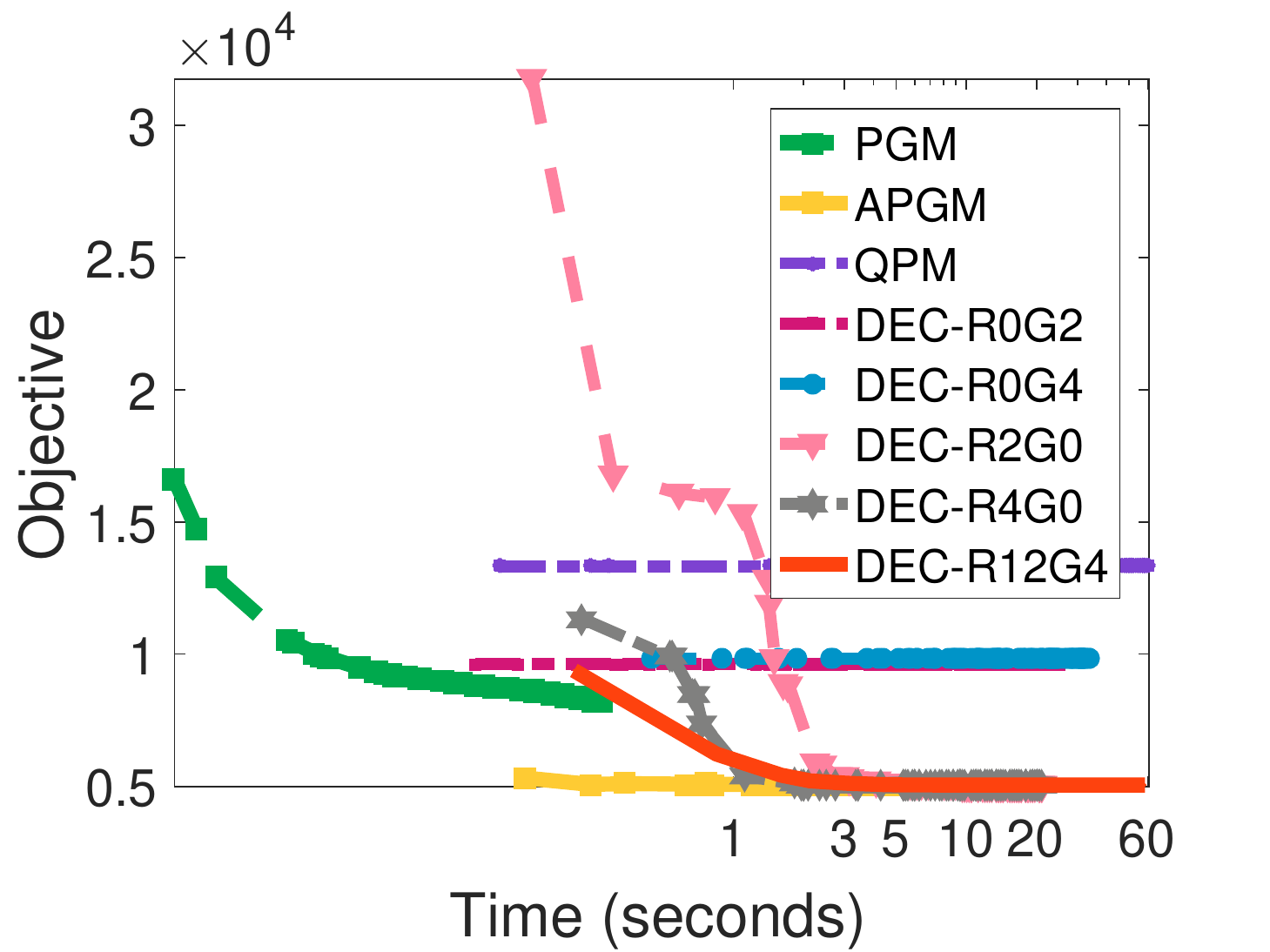}\vspace{-6pt} \caption{\scriptsize s=20 on e2006-5000-1024 }\end{subfigure}\ghs
      \begin{subfigure}{\fourfigwid}\includegraphics[width=\objimgwid,height=\objimghei]{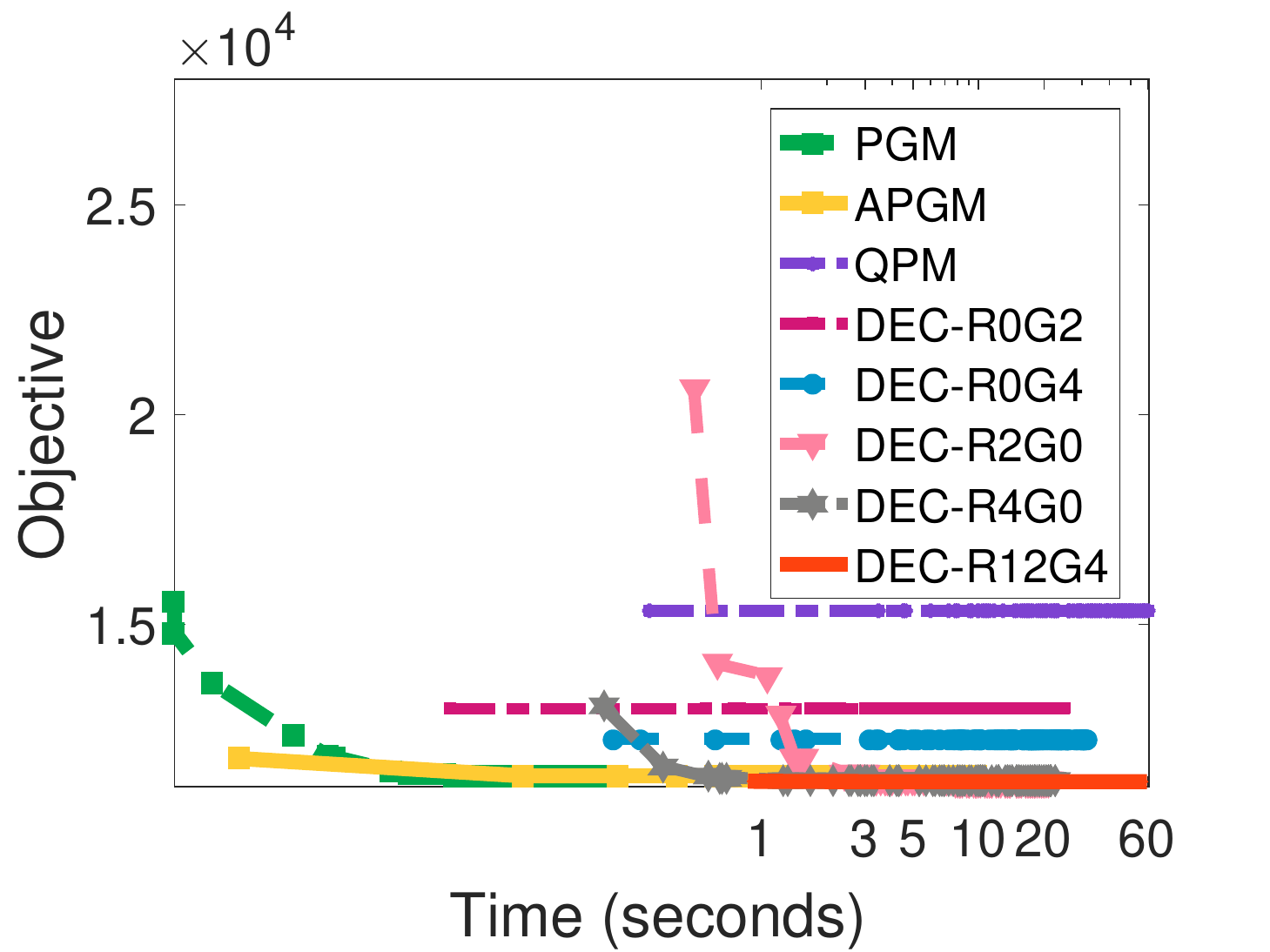}\vspace{-6pt} \caption{\scriptsize s=40 on e2006-5000-1024 }\end{subfigure}\\

     \begin{subfigure}{\fourfigwid}\includegraphics[width=\objimgwid,height=\objimghei]{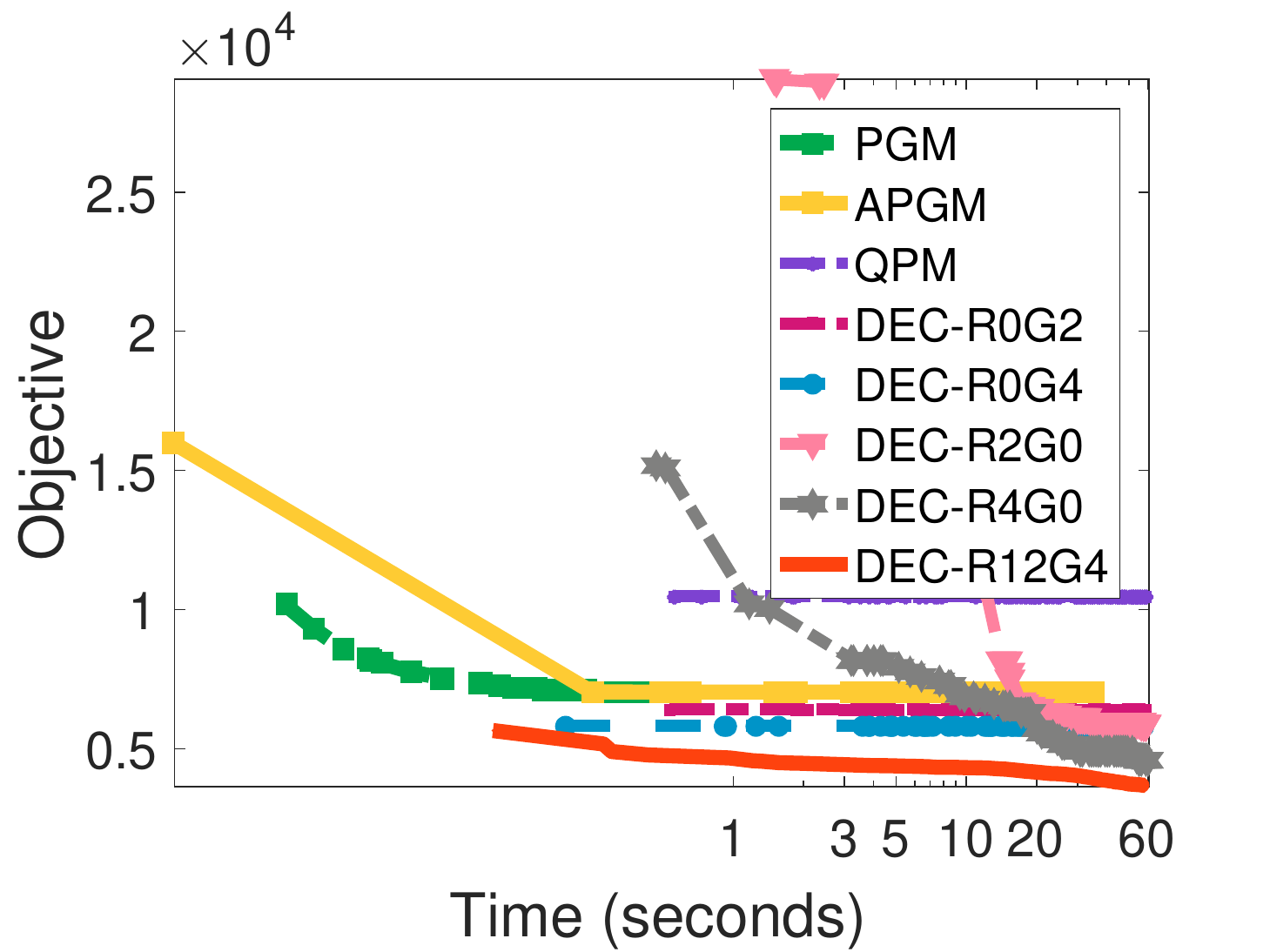}\vspace{-6pt} \caption{\scriptsize s=20 on e2006-5000-2048 }\end{subfigure}\ghs
     \begin{subfigure}{\fourfigwid}\includegraphics[width=\objimgwid,height=\objimghei]{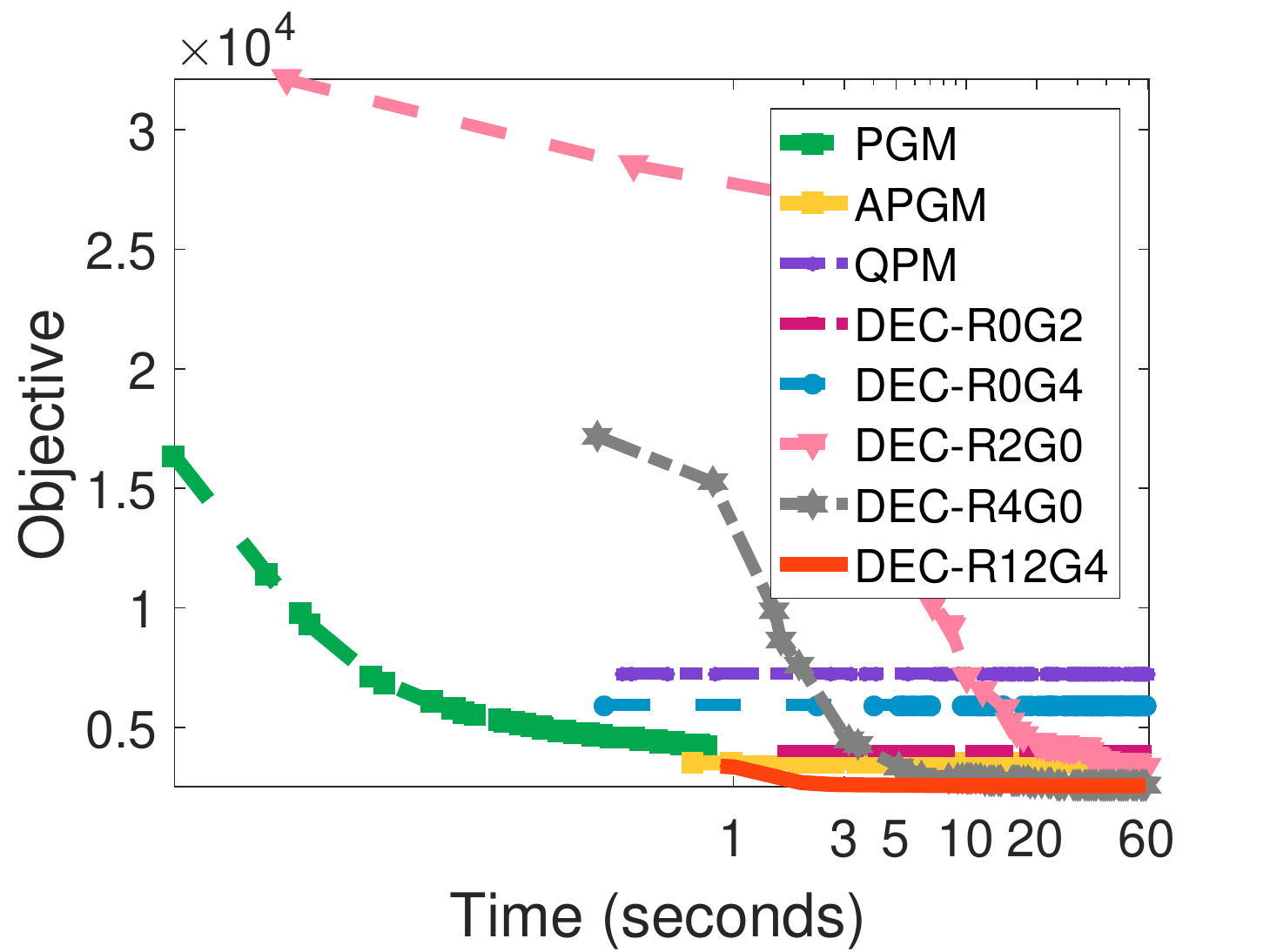}\vspace{-6pt} \caption{ \scriptsize s=40 on e2006-5000-2048} \end{subfigure}\ghs
      \begin{subfigure}{\fourfigwid}\includegraphics[width=\objimgwid,height=\objimghei]{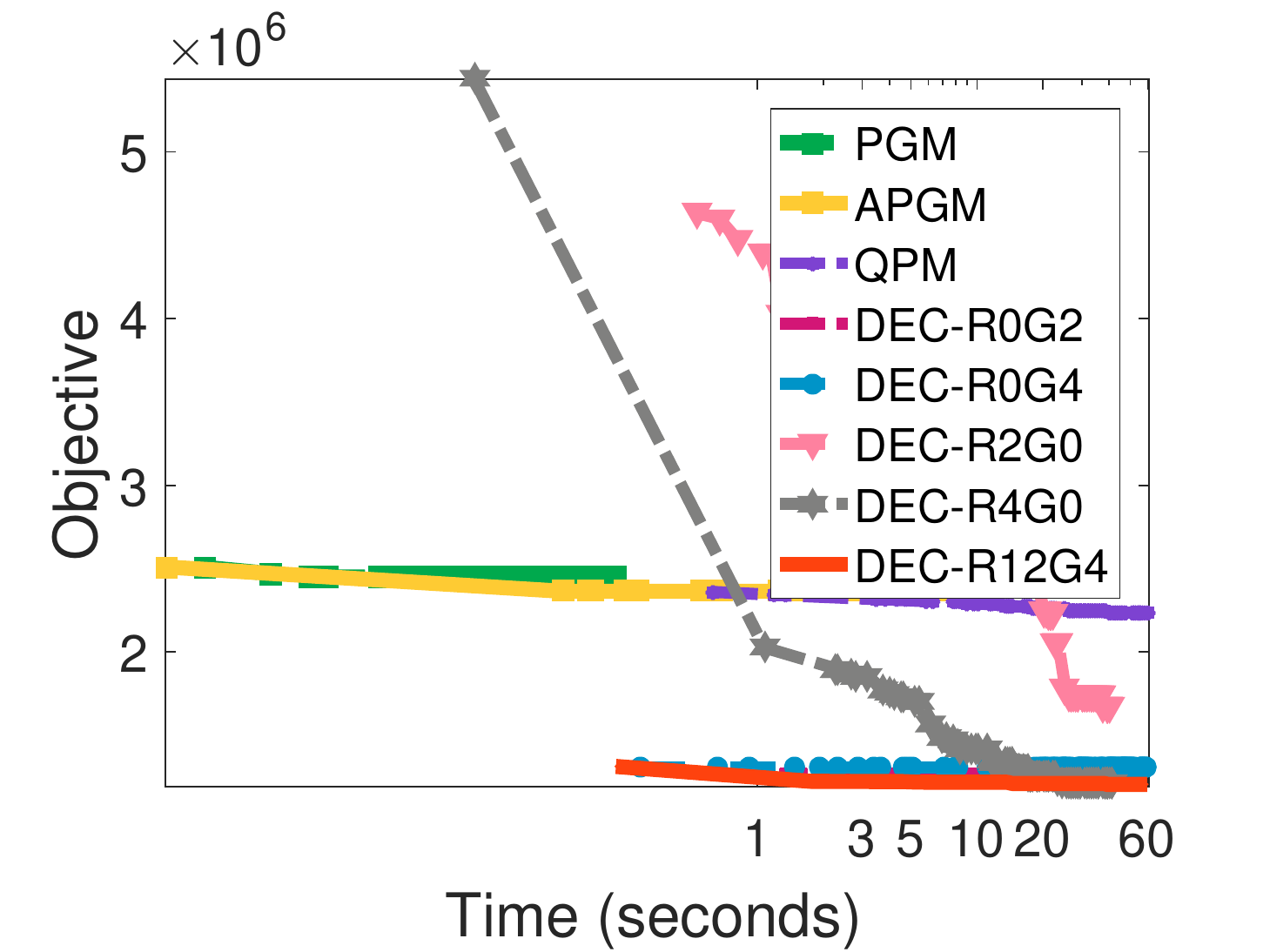}\vspace{-6pt} \caption{\scriptsize s=20 on random-256-1024-C }\end{subfigure}\ghs
      \begin{subfigure}{\fourfigwid}\includegraphics[width=\objimgwid,height=\objimghei]{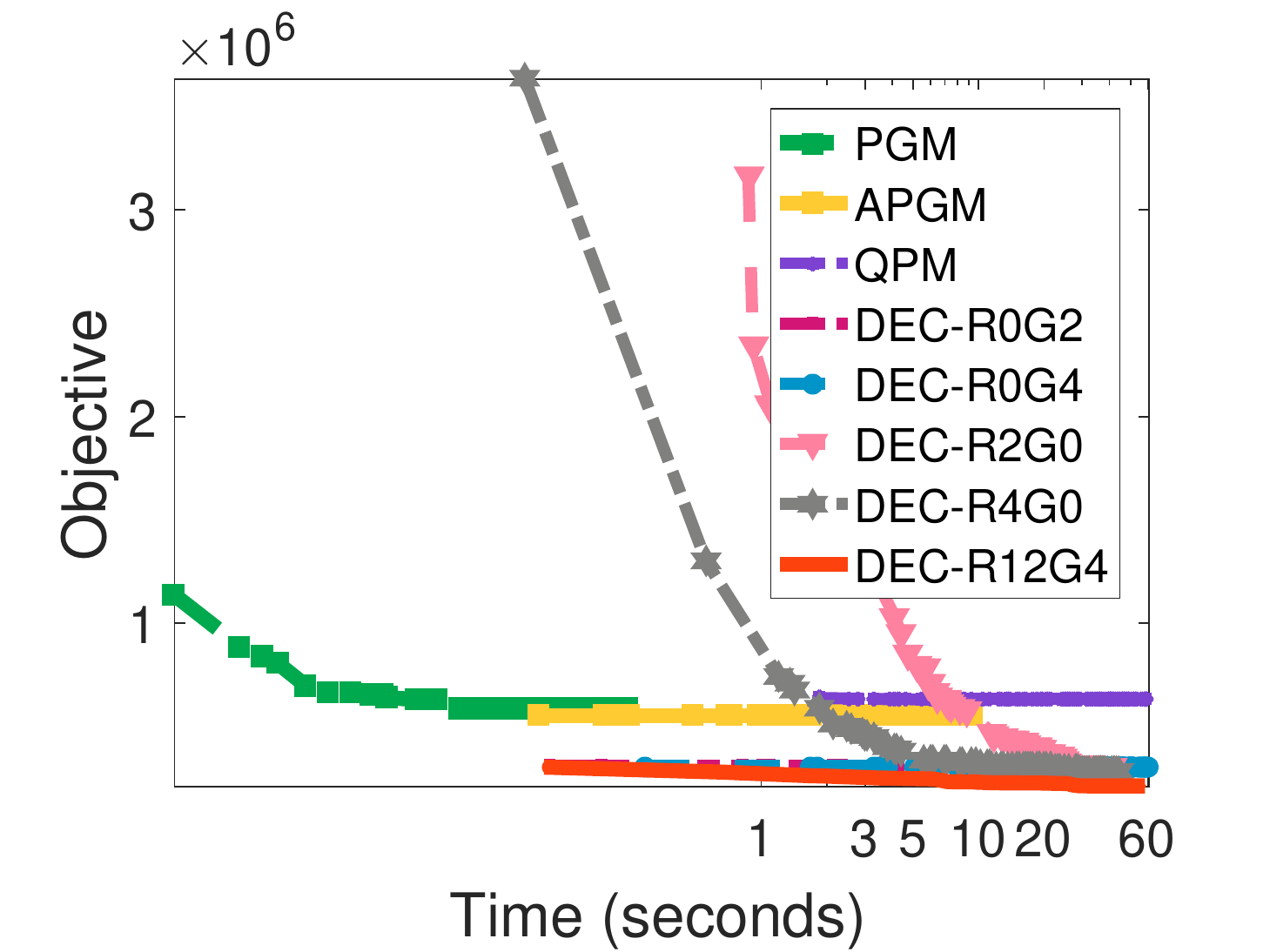}\vspace{-6pt} \caption{\scriptsize s=40 on random-256-1024-C }\end{subfigure}\\

     \begin{subfigure}{\fourfigwid}\includegraphics[width=\objimgwid,height=\objimghei]{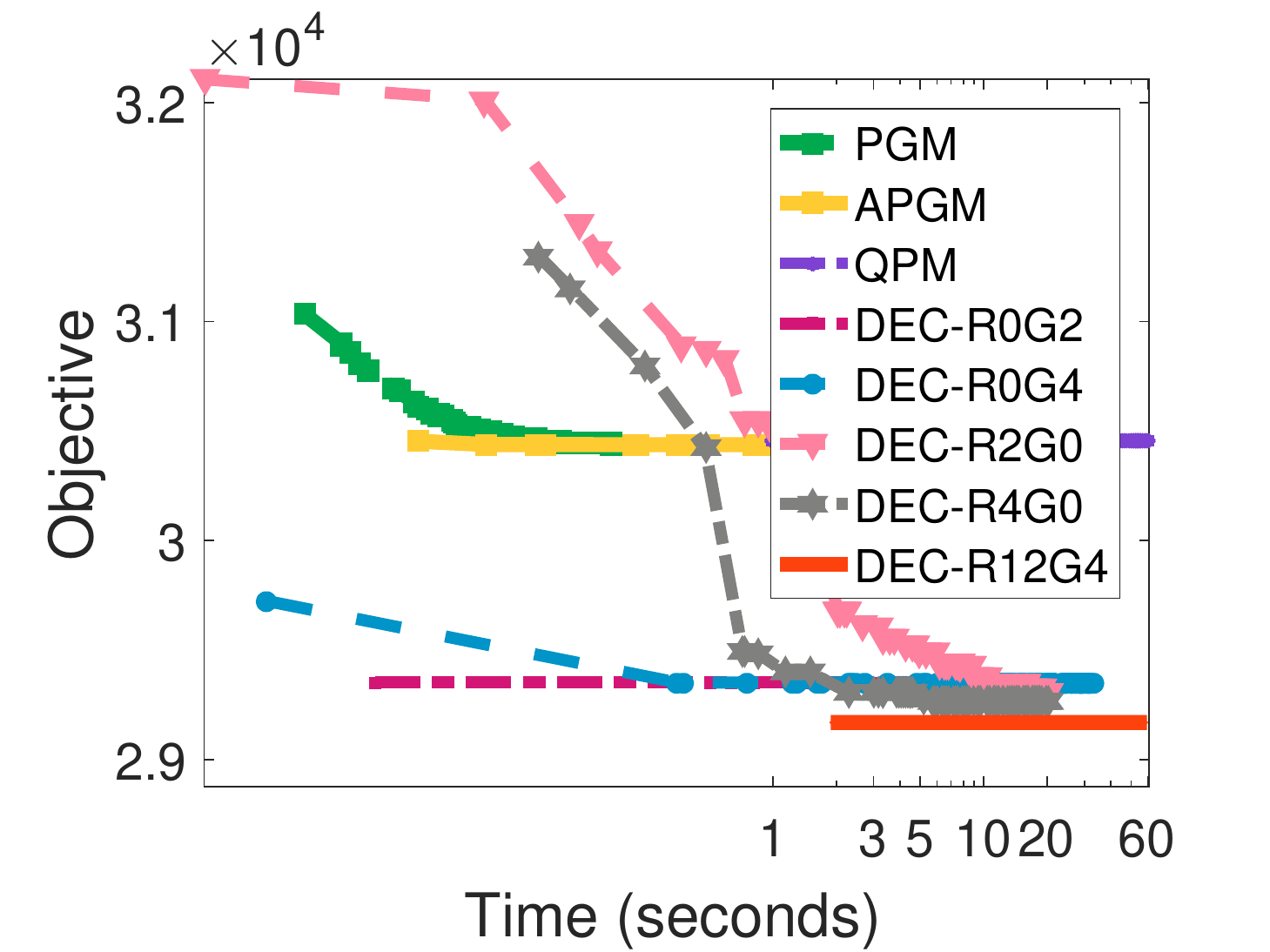}\vspace{-6pt} \caption{\scriptsize s=20 on e2006-5000-1024-C}\end{subfigure}\ghs
     \begin{subfigure}{\fourfigwid}\includegraphics[width=\objimgwid,height=\objimghei]{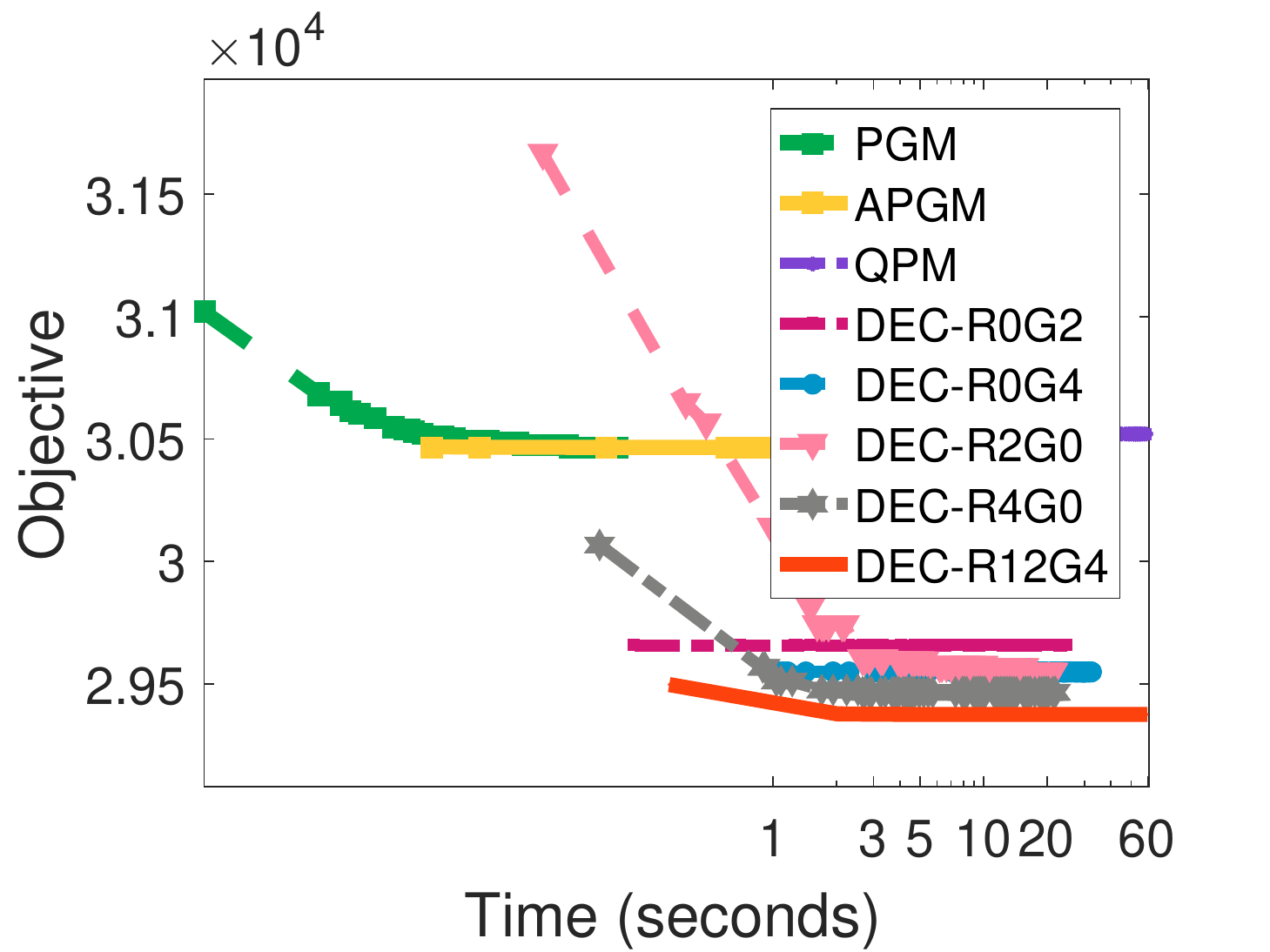}\vspace{-6pt} \caption{ \scriptsize s=40 on e2006-5000-1024-C } \end{subfigure}\ghs
      \begin{subfigure}{\fourfigwid}\includegraphics[width=\objimgwid,height=\objimghei]{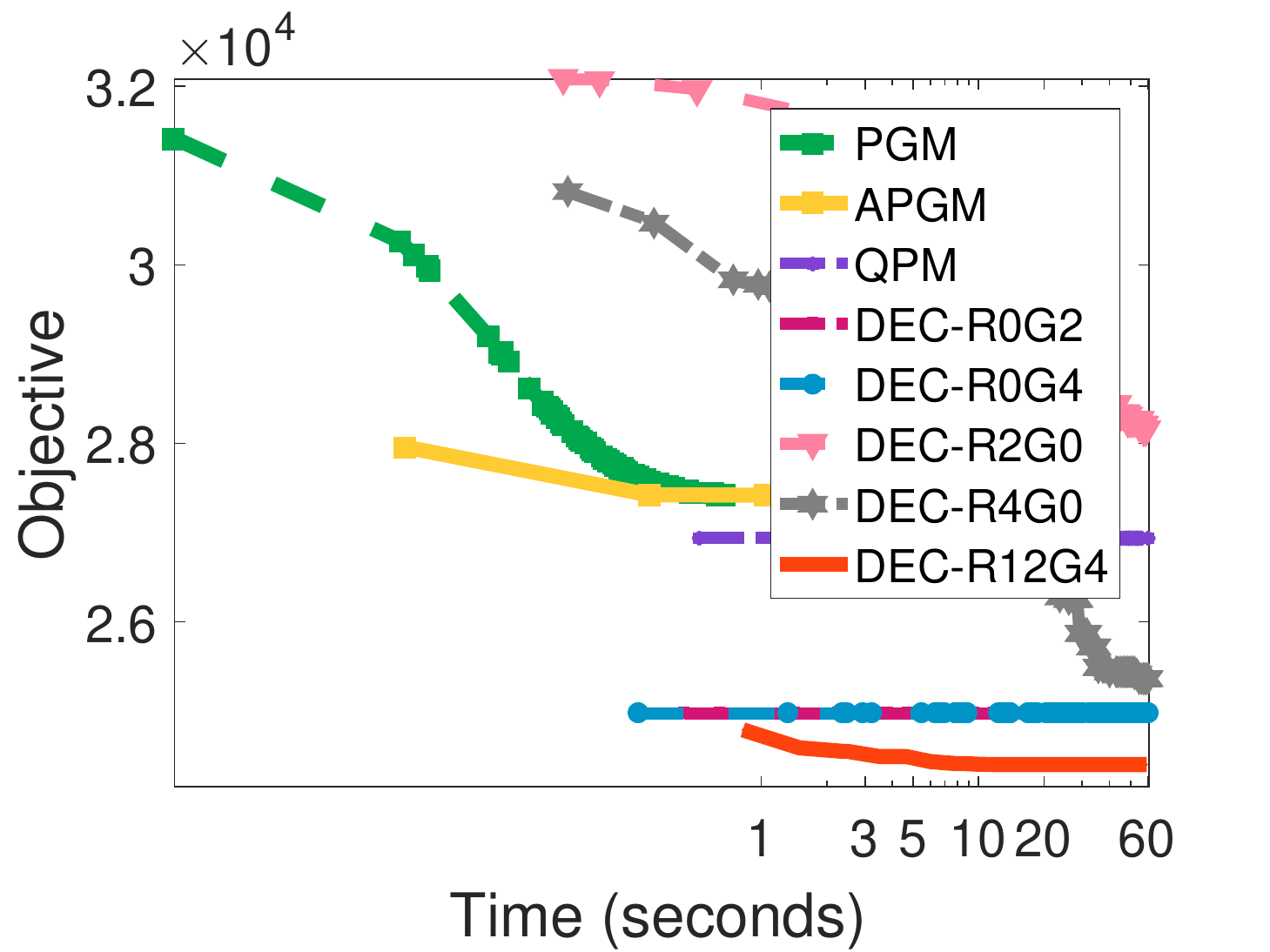}\vspace{-6pt} \caption{\scriptsize s=20 on e2006-5000-2048-C}\end{subfigure}\ghs
      \begin{subfigure}{\fourfigwid}\includegraphics[width=\objimgwid,height=\objimghei]{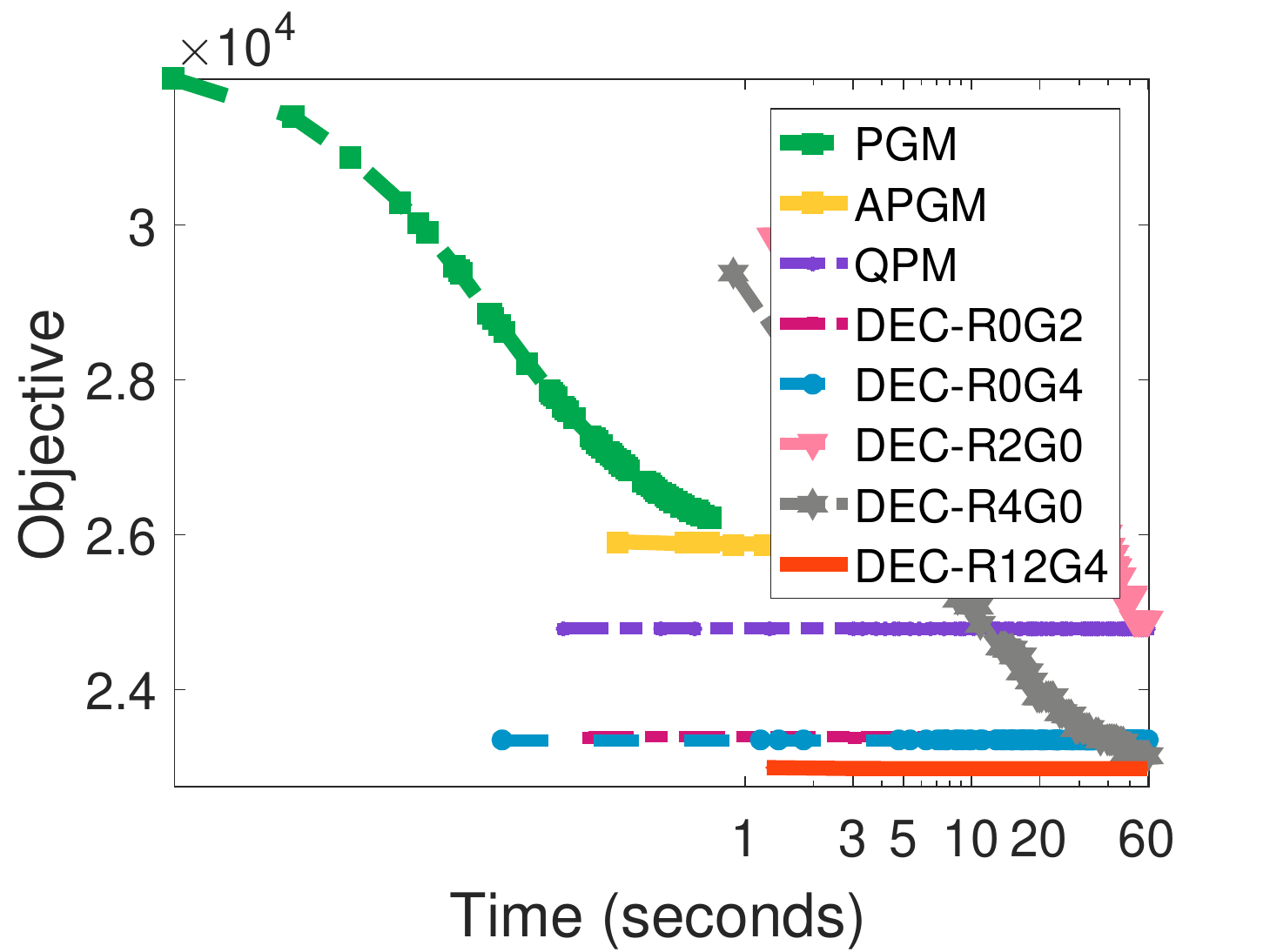}\vspace{-6pt} \caption{\scriptsize s=40 on e2006-5000-2048-C }\end{subfigure}\\

\centering
\caption{Convergence curve and computional efficiency for solving sparsity constrained least squares problems on different data sets with different $s$. }
\label{fig:1}
\end{figure*}

\begin{figure*} [!t]
\centering
      \begin{subfigure}{\fourfigwid}\includegraphics[width=\objimgwid,height=\objimghei]{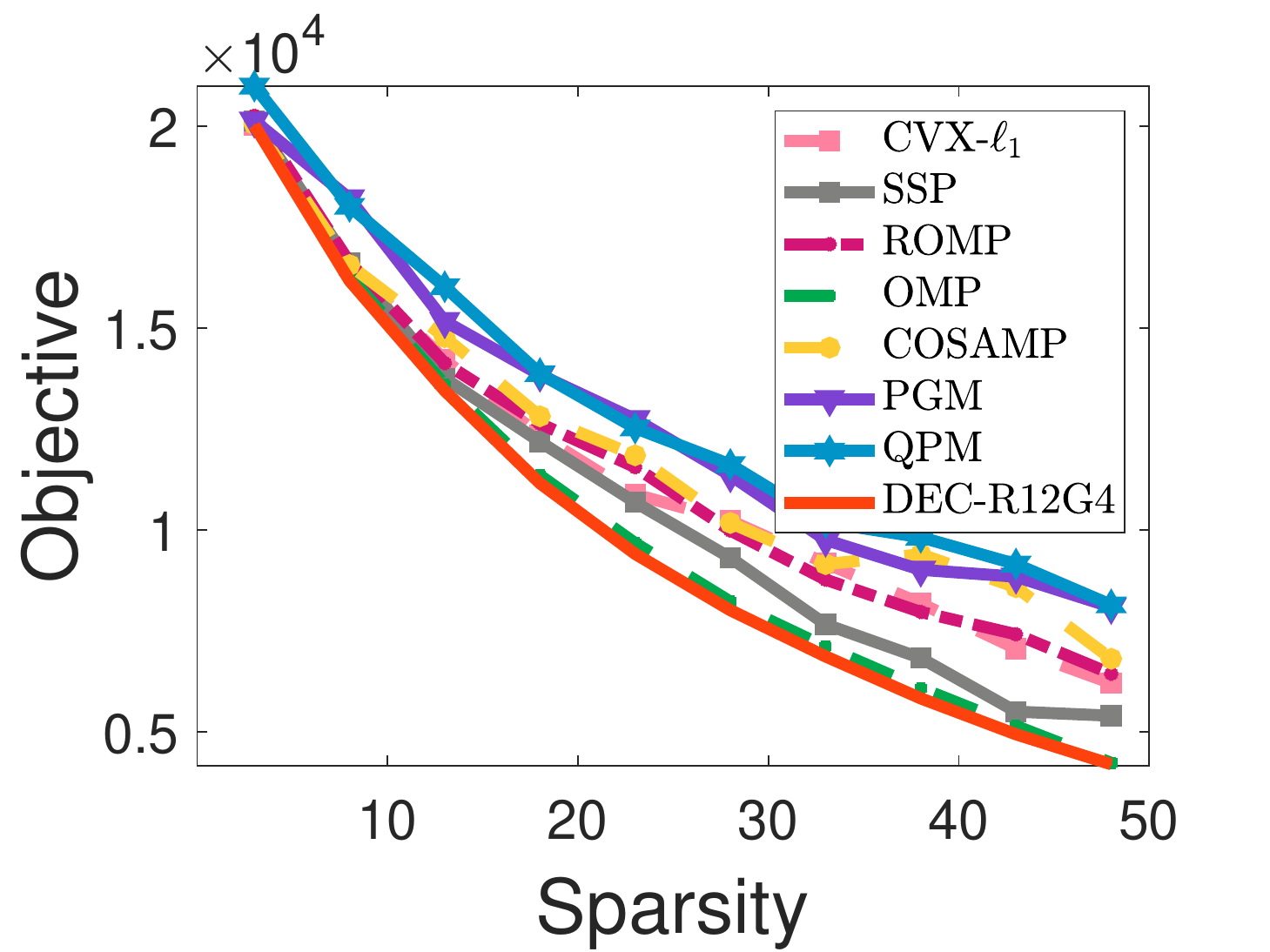}\vspace{-6pt} \caption{\scriptsize random-256-1024 }\end{subfigure}\ghs
     \begin{subfigure}{\fourfigwid}\includegraphics[width=\objimgwid,height=\objimghei]{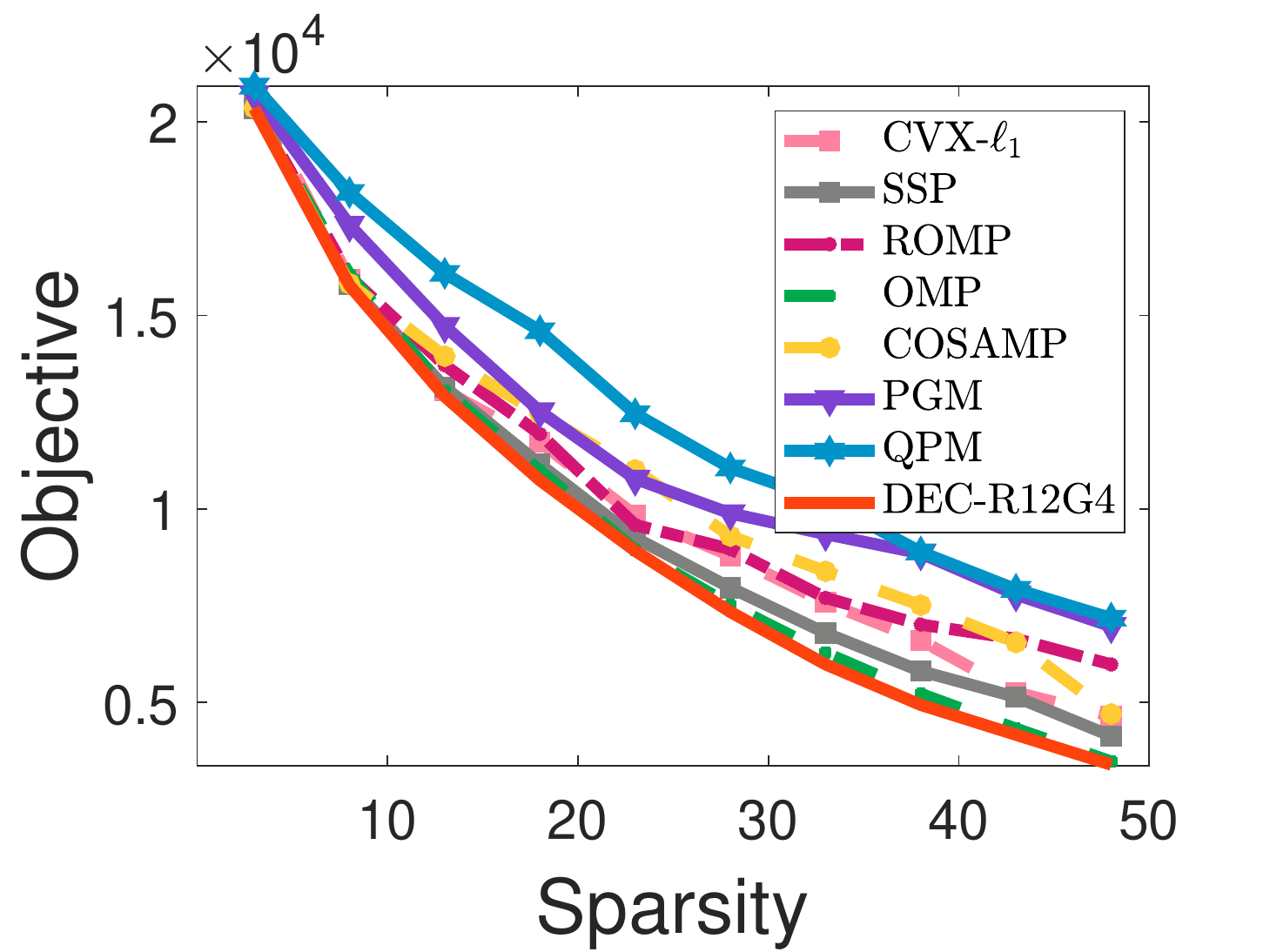}\vspace{-6pt} \caption{ \scriptsize random-256-2048} \end{subfigure}\ghs
      \begin{subfigure}{\fourfigwid}\includegraphics[width=\objimgwid,height=\objimghei]{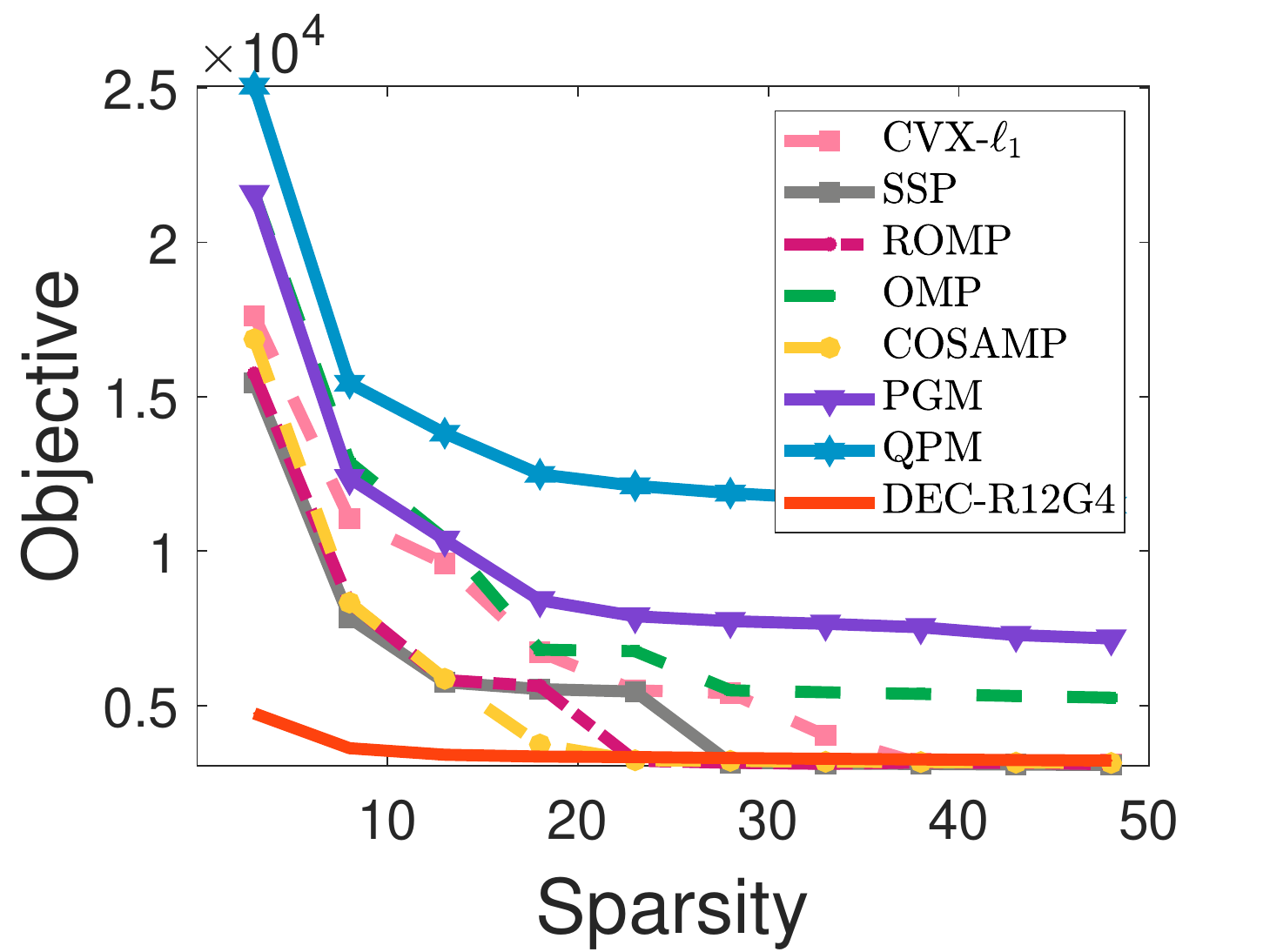}\vspace{-6pt} \caption{\scriptsize e2006-5000-1024 }\end{subfigure}\ghs
      \begin{subfigure}{\fourfigwid}\includegraphics[width=\objimgwid,height=\objimghei]{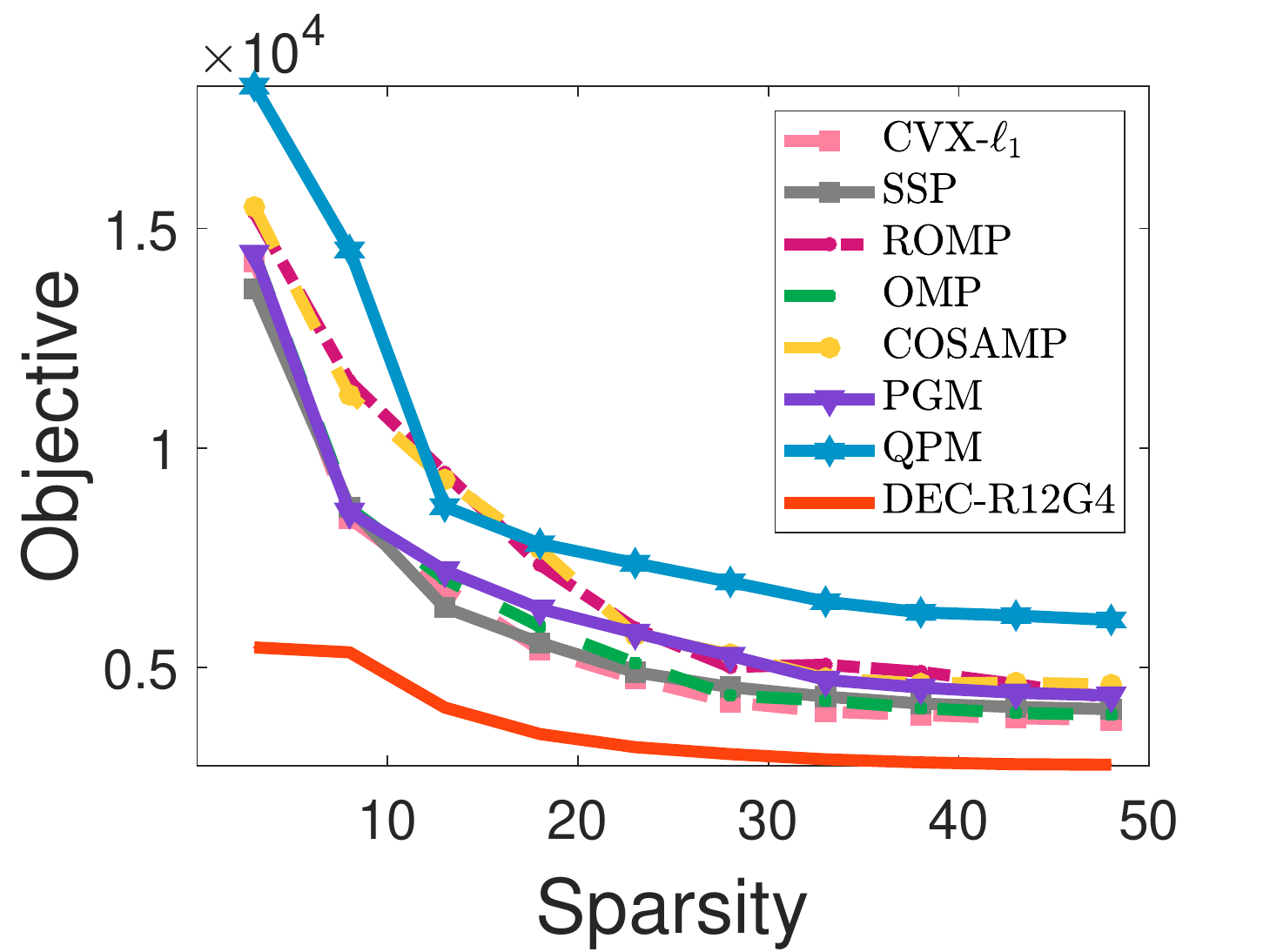}\vspace{-6pt} \caption{\scriptsize e2006-5000-2048}\end{subfigure}\\

      \begin{subfigure}{\fourfigwid}\includegraphics[width=\objimgwid,height=\objimghei]{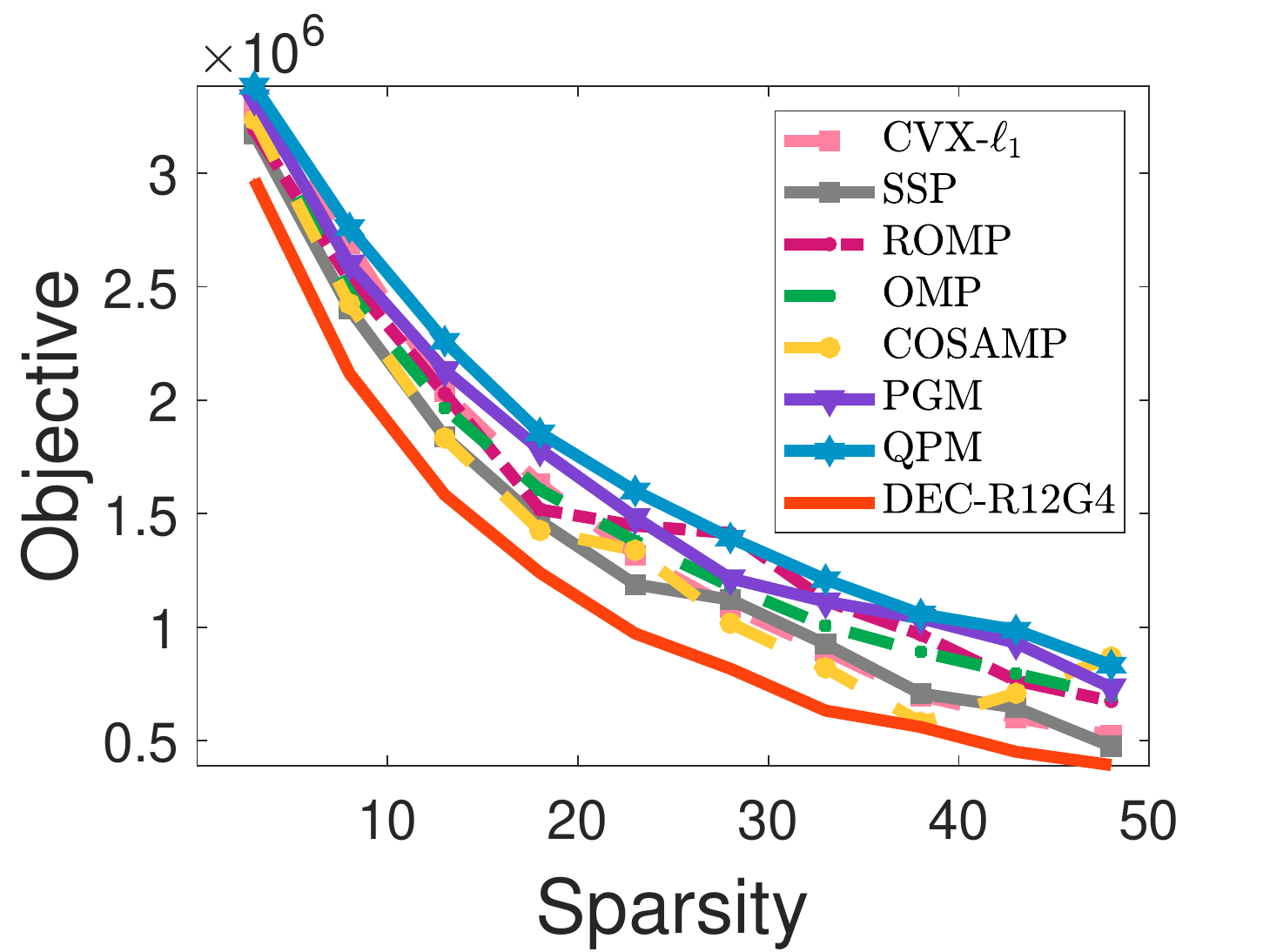}\vspace{-6pt} \caption{ \scriptsize random-256-1024-C }\end{subfigure}\ghs
       \begin{subfigure}{\fourfigwid}\includegraphics[width=\objimgwid,height=\objimghei]{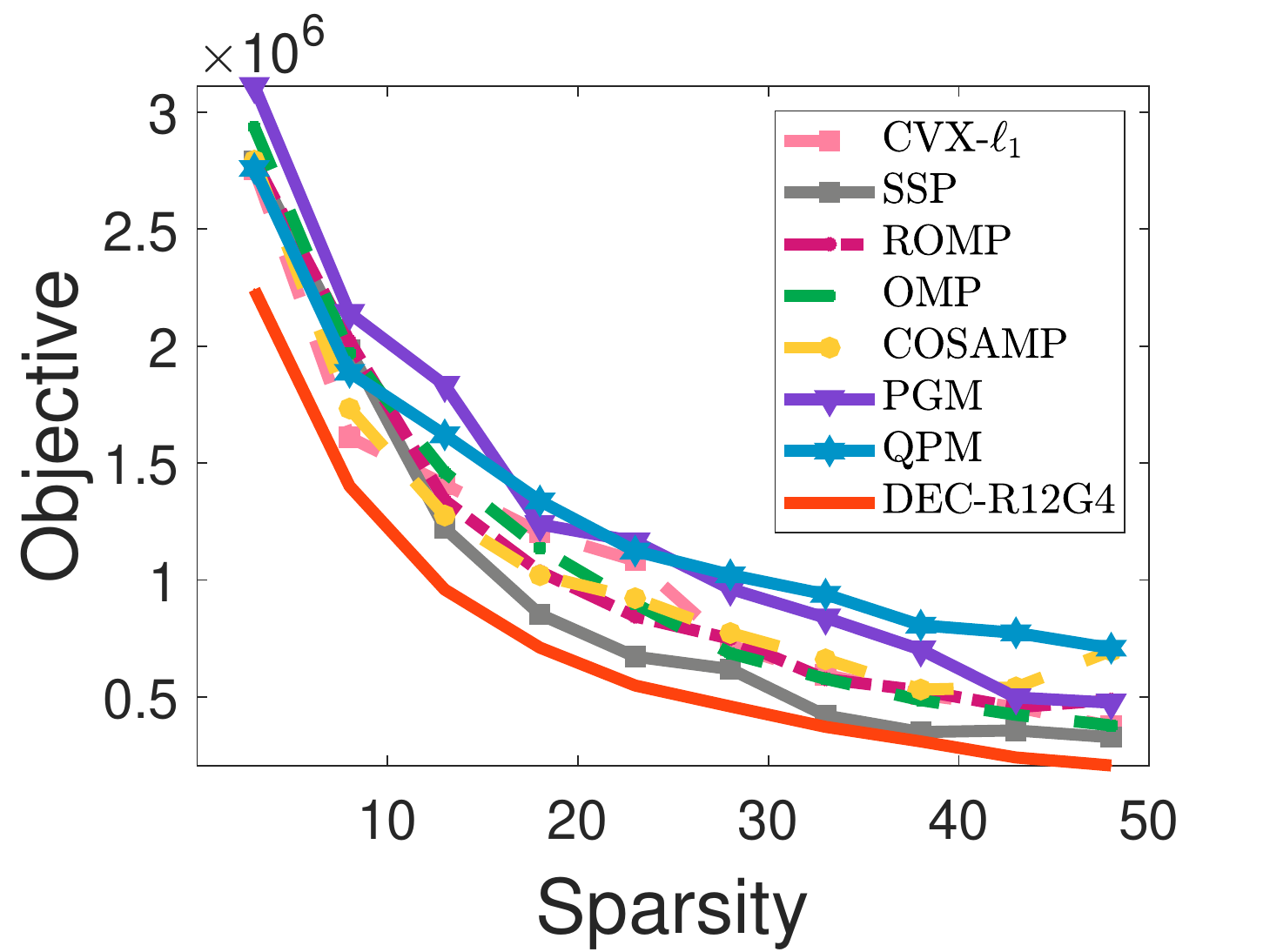}\vspace{-6pt} \caption{ \scriptsize random-256-2048-C}\end{subfigure}\ghs
      \begin{subfigure}{\fourfigwid}\includegraphics[width=\objimgwid,height=\objimghei]{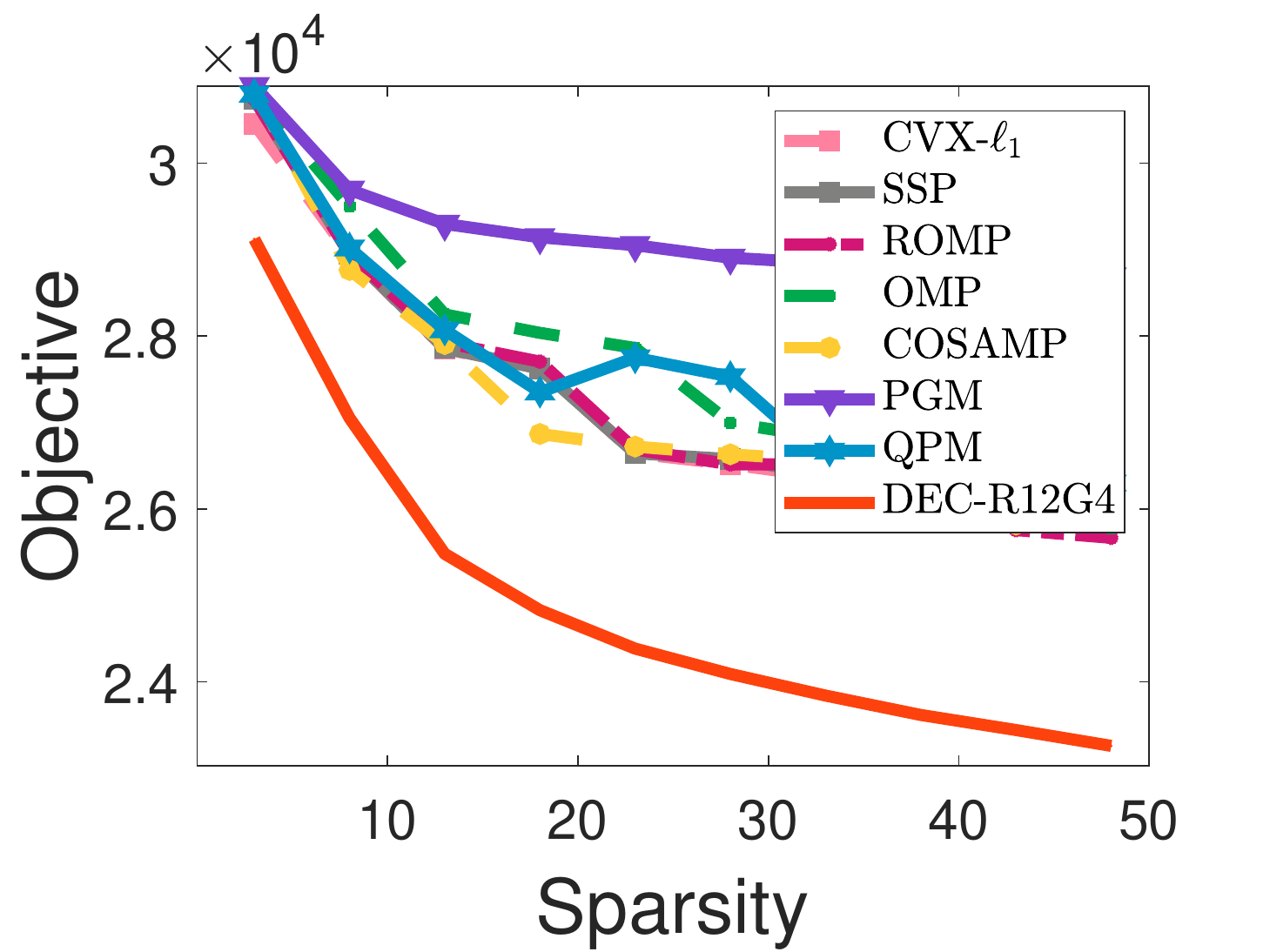}\vspace{-6pt} \caption{ \scriptsize e2006-5000-1024-C}\end{subfigure}\ghs
      \begin{subfigure}{\fourfigwid}\includegraphics[width=\objimgwid,height=\objimghei]{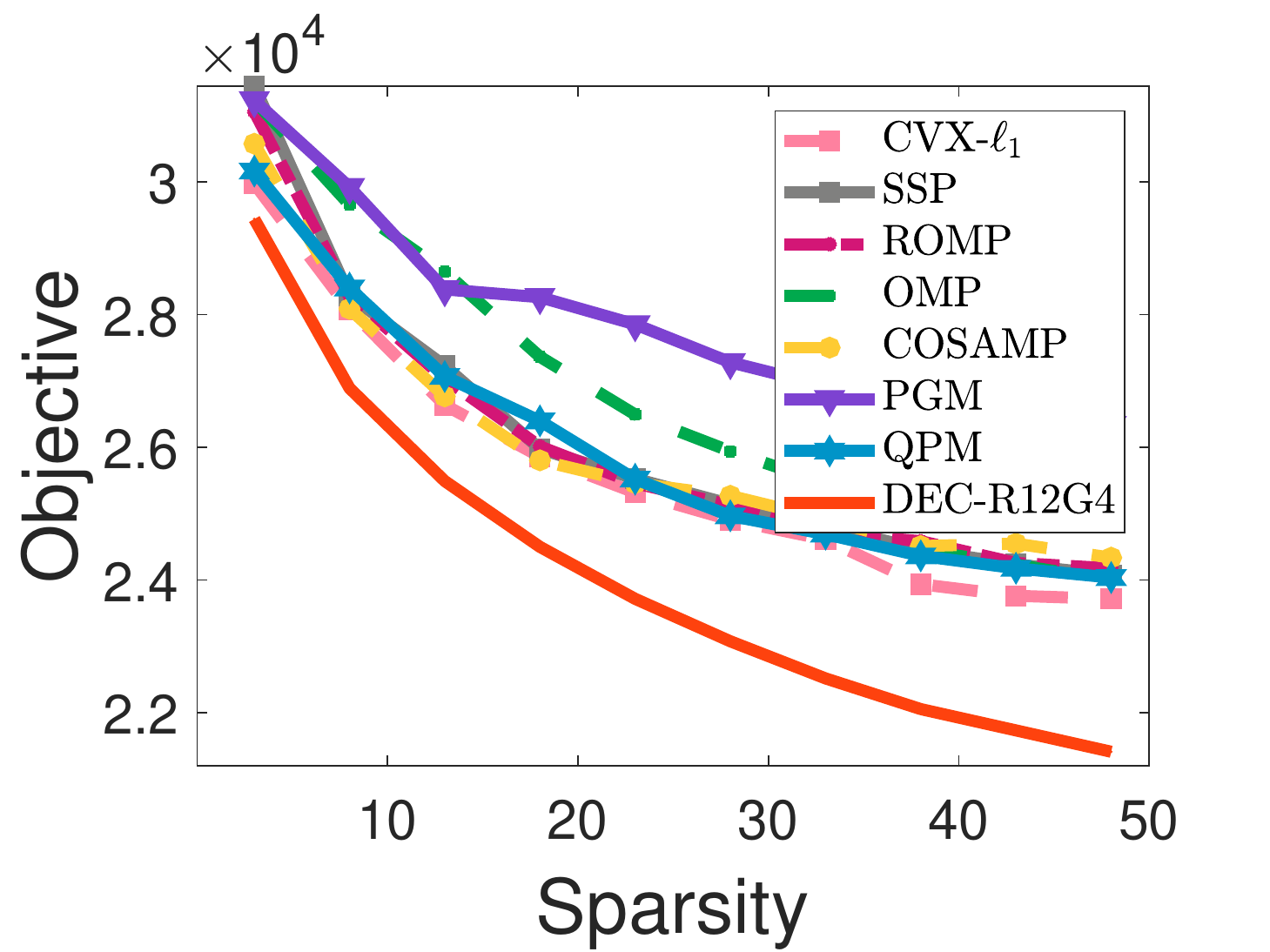}\vspace{-6pt} \caption{ \scriptsize e2006-5000-2048-C}\end{subfigure}\\

\centering
\caption{Experimental results on sparsity constrained least squares problems on different data sets with varying the sparsity of the solution. }
\label{fig2}
\vspace{-3pt}
\end{figure*}

\section{Discussions} \label{sect:dis}

This section provides additional discussions for the proposed method.

  \bbb{$\blacktriangleright$ When the objective function is complicated.} In step (S2) of the proposed algorithm, a global solution is to be found for the subproblem. When $f(\cdot)$ is simple (e.g., a quadratic function), we can find efficient and exact solutions to the subproblems. We now consider the situation when $f$ is complicated (e.g., logistic regression, maximum entropy models). One can still find a quadratic majorizer $Q(\bbb{x},\bbb{z})$ for the convex smooth function $f(\bbb{x})$ with
\beq
\forall~\bbb{z},~\bbb{x},~f(\bbb{x}) \leq Q(\bbb{x},\bbb{z}) \triangleq f(\bbb{z}) + (\bbb{x}-\bbb{z})^{\intercal}\nabla f(\bbb{z})+\nn\\
~~~~~~\tfrac{1}{2}(\bbb{x}-\bbb{z})^{\intercal}\bbb{M}(\bbb{z})(\bbb{x}-\bbb{z}),~~\bbb{M}(\bbb{z})\succ \nabla f^2(\bbb{z}). \nn
\eeq
\noi By minimizing the upper bound of $f(\bbb{x})$ (i.e., the quadratic surrogate function) at the current estimate $\bbb{x}^t$, i.e.,
\beq
\bbb{x}^{t+1}\Leftarrow \arg\min_{\bbb{x}}~Q(\bbb{x},\bbb{x}^{t})+h(\bbb{x}),\nn
\eeq
 \noi we can drive the objective downward until a stationary point is reached. We will obtain a stationary point $\ddot{\bbb{x}}$ satisfying:
 \beq
 \ddot{\bbb{x}} = \arg \min_{\bbb{z}}~h(\bbb{z}) + f(\ddot{\bbb{x}}) + (\bbb{z}-\ddot{\bbb{x}})^{\intercal}\nabla f(\ddot{\bbb{x}})  \nn\\
 + \fractt{1}{2} (\bbb{z}-\ddot{\bbb{x}})^{\intercal}  \bbb{M}(\ddot{\bbb{x}}) (\bbb{z}-\ddot{\bbb{x}}),~s.t.~\ddot{\bbb{x}}_{\bar{B}}=\bbb{z}_{\bar{B}}\nn
 \eeq
  \noi for all $\bar{B}$. We denote this optimality condition as the \underline{Newton} block-$k$ stationary point. It is weaker than the \underline{full} block-$k$ stationary point as in Definition \ref{def:block:k}. However, it is stronger than the L-Stationary Point as in Definition \ref{def:L:stationary}.


  \bbb{$\blacktriangleright$ Computational efficiency.} Block coordinate descent is shown to be very efficient for solving convex problems (e.g., support vector machines \cite{chang2008coordinate,hsieh2008dual}, LASSO problems \cite{tseng2009coordinate}, nonnegative matrix factorization \cite{hsieh2011fast}). The main difference of our block coordinate descent from existing ones is that our method needs to solve a small-sized NP-hard subproblem globally which takes subexponential
time $\mathcal{O}(2^k)$. As a result, our algorithm finds a block-$k$ approximation solution
for the original NP-hard problem within $\mathcal{O}(2^k)$ time. When $k$ is large, it is hard to enumerate the full binary tree since the subproblem is equally NP-hard. However, $k$ is relatively small in practice (e.g., 2 to 20). In addition, real-world applications often have some special (e.g., unbalanced, sparse, local) structure and block-$k$ stationary point could also be the global stationary point (refer to Table \ref{tab:optimality} of this paper).



\section{Experimental Validation} \label{sect:exp}
This section demonstrates the effectiveness of our algorithm on two sparse optimization tasks, namely the sparse regularized least squares problem and the sparsity constrained least squares problem.

Given a design matrix $\bbb{A} \in \mathbb{R}^{m\times n}$ and an observation vector $\bbb{b} \in \mathbb{R}^m$, we solve the following optimization problem:
\beq
\begin{split}
\textstyle
 \min_{\bbb{x}}~\fractt{1}{2}\|\bbb{Ax}-\bbb{b}\|_2^2,~s.t.~\|\bbb{x}\|_0 \leq s, ~~~\nn\\
\textstyle \text{or}~~~\min_{\bbb{x}}\fractt{1}{2}\|\bbb{Ax}-\bbb{b}\|_2^2 +  \lambda\|\bbb{x}\|_0,
\end{split}
\eeq

\noi where  $s$ and $\lambda$ are given parameters.


\bbb{Experimental Settings.} We use \textbf{DEC-R$i$G$j$} to denote our block decomposition method along with selecting $i$ coordinates using the \bbb{R}andom strategy and $j$ coordinates using the \bbb{G}reedy strategy. We keep a record of the relative changes of the objective function values by $r_t = (f(\bbb{x}^t)-f(\bbb{x}^{t+1}))/f(\bbb{x}^t)$. We let \textbf{DEC} run up to $T$ iterations and stop it at iteration $t<T$ if $\text{mean}([{r}_{t-\text{min}(t,\varrho)+1},~{r}_{t-min(t,\varrho)+2},...,{r}_t]) \leq \epsilon$. We use the default value $(\theta,~\epsilon,~\varrho,~T)=(10^{-3},~10^{-5},~50,~1000)$ for \textbf{DEC}. All codes were implemented in Matlab on an Intel 3.20GHz CPU with 8 GB RAM. We measure the quality of the solution by comparing the objective values for different methods. Note that although
\textbf{DEC} uses the randomized strategy to find the working set, we can always measure the quality of the solution by computing the deterministic objective value.

\begin{table}[!h]
\scalebox{0.94}{\begin{tabular}{|c|c|c|c|c|c|}
\hline
&{\tiny  PGM-$\ell_0$} & {\tiny APGM-$\ell_0$} &  {\tiny PGM-$\ell_1$} &{\tiny PGM-$\ell_p$} &{\tiny DEC-R10G2} \\
\hline
  \hline
  \multicolumn{6}{|>{\columncolor{mycyana}}c|}{\centering results on random-256-1024 } \\
$\lambda = 10^0$ &\cthree{6.9e+2} & 2.4e+4& 7.8e+2& \cone{4.0e+2} & \ctwo{4.8e+2} \\
 $\lambda = 10^1$ &\cthree{2.3e+3} & 3.8e+4& 3.3e+3& \cone{1.9e+3} & \ctwo{2.2e+3} \\
 $\lambda = 10^2$ &2.0e+4& 1.3e+5& \cthree{1.8e+4} & \ctwo{1.1e+4} & \cone{9.4e+3} \\
 $\lambda = 10^3$ &2.5e+4& 1.0e+6& \ctwo{2.4e+4} & \cthree{2.4e+4} & \cone{2.4e+4} \\
\hline\hline
  \multicolumn{6}{|>{\columncolor{mycyanb}}c|}{\centering results on random-256-2048 } \\
\hline
$\lambda = 10^0$ &\cthree{1.3e+3} & 2.7e+4& 1.4e+3& \ctwo{6.0e+2} & \cone{5.4e+2} \\
 $\lambda = 10^1$ &\cthree{2.9e+3} & 4.5e+4& 4.9e+3& \cone{2.2e+3} & \ctwo{2.2e+3} \\
 $\lambda = 10^2$ &2.2e+4& 2.3e+5& \cthree{2.1e+4} & \ctwo{1.1e+4} & \cone{9.5e+3} \\
 $\lambda = 10^3$ &\cthree{2.7e+4} & 2.1e+6& \ctwo{2.6e+4} & \cthree{2.7e+4} & \cone{2.6e+4} \\
\hline\hline
  \multicolumn{6}{|>{\columncolor{mycyanc}}c|}{\centering results on e2006-5000-1024 } \\
\hline
 $\lambda = 10^0$ &\ctwo{8.5e+3} & 3.3e+4& \cthree{1.1e+4} & 1.8e+4& \cone{7.3e+3} \\
 $\lambda = 10^1$ &\ctwo{9.4e+3} & 4.2e+4& \cthree{3.2e+4} & \cthree{3.2e+4} & \cone{8.6e+3} \\
 $\lambda = 10^2$ &\ctwo{3.2e+4} & \cthree{1.3e+5} & \ctwo{3.2e+4} & \ctwo{3.2e+4} & \cone{1.3e+4} \\
 $\lambda = 10^3$ &\ctwo{1.8e+4} & 1.1e+6& \cthree{3.2e+4} & \cthree{3.2e+4} & \cone{1.1e+4} \\
\hline\hline
  \multicolumn{6}{|>{\columncolor{mycyand}}c|}{\centering results on e2006-5000-2048 } \\
\hline
 $\lambda = 10^0$ &\ctwo{3.1e+3} & 3.4e+4& \cthree{4.4e+3} & 1.4e+4& \cone{2.6e+3} \\
 $\lambda = 10^1$ &\ctwo{5.2e+3} & 5.3e+4& \cthree{1.2e+4} & 1.2e+4& \cone{4.5e+3} \\
 $\lambda = 10^2$ &\ctwo{3.2e+4} & \cthree{2.4e+5} & \ctwo{3.2e+4} & \ctwo{3.2e+4} & \cone{7.0e+3} \\
 $\lambda = 10^3$ &\ctwo{1.8e+4} & 2.1e+6& \cthree{3.2e+4} & \cthree{3.2e+4} & \cone{1.3e+4} \\
\hline\hline
  \multicolumn{6}{|>{\columncolor{mycyana}}c|}{\centering results on random-256-1024-C } \\
\hline
 $\lambda = 10^0$ &\ctwo{9.6e+2} & 5.7e+6& 1.0e+3& \cthree{1.0e+3} & \cone{8.9e+2} \\
 $\lambda = 10^1$ &\ctwo{8.1e+3} & 3.5e+6& 1.0e+4& \cthree{8.2e+3} & \cone{7.3e+3} \\
 $\lambda = 10^2$ &\cthree{5.8e+4} & 6.2e+6& 8.9e+4& \ctwo{5.4e+4}  & \cone{5.1e+4} \\
 $\lambda = 10^3$ &\cthree{2.5e+5} & 5.3e+6& 3.7e+5& \ctwo{2.2e+5}  &  \cone{2.0e+5}\\
\hline\hline
  \multicolumn{6}{|>{\columncolor{mycyanb}}c|}{\centering results on random-256-2048-C } \\
\hline
 $\lambda = 10^0$ &\ctwo{1.9e+3} & 5.7e+6& 2.0e+3& \cthree{1.9e+3} & \cone{1.2e+3} \\
 $\lambda = 10^1$ &\cthree{1.7e+4} & 7.7e+6& 2.0e+4& \ctwo{1.6e+4} & \cone{9.2e+3} \\
 $\lambda = 10^2$ &\cthree{8.4e+4} & 4.2e+6& 1.6e+5& \ctwo{6.4e+4} & \cone{5.3e+4} \\
 $\lambda = 10^3$ & \cthree{2.5e+5}& 9.6e+6& 6.3e+5& \ctwo{2.5e+5} & \cone{2.4e+5}  \\
\hline\hline
  \multicolumn{6}{|>{\columncolor{mycyanc}}c|}{\centering results on e2006-5000-1024-C } \\
\hline
 $\lambda = 10^0$ &3.0e+4& 3.3e+4& \ctwo{2.8e+4} & \cthree{2.9e+4} & \cone{2.2e+4} \\
 $\lambda = 10^1$ &\ctwo{3.2e+4} & \cthree{4.2e+4} & \ctwo{3.2e+4} & \ctwo{3.2e+4} & \cone{2.3e+4} \\
 $\lambda = 10^2$ &\cthree{3.2e+4} & 1.3e+5& \ctwo{3.2e+4} & 3.2e+4& \cone{2.9e+4} \\
 $\lambda = 10^3$ &\cone{3.2e+4} & \ctwo{1.1e+6} & \cone{3.2e+4} & \cone{3.2e+4} & \cone{3.2e+4} \\
\hline\hline
  \multicolumn{6}{|>{\columncolor{mycyand}}c|}{\centering results on e2006-5000-2048-C } \\
\hline
 $\lambda = 10^0$ &2.9e+4& 3.4e+4& \ctwo{2.6e+4} & \cthree{2.7e+4} & \cone{1.7e+4} \\
 $\lambda = 10^1$ &\ctwo{3.2e+4} & 5.3e+4& \cthree{3.2e+4} & \ctwo{3.2e+4} & \cone{2.1e+4} \\
 $\lambda = 10^2$ &\ctwo{3.2e+4} & \cthree{2.4e+5} & \ctwo{3.2e+4} & \ctwo{3.2e+4} & \cone{2.7e+4} \\
 $\lambda = 10^3$ &\cone{3.2e+4} & \ctwo{2.1e+6} & \cone{3.2e+4} & \cone{3.2e+4} & \cone{3.2e+4} \\
 \hline
\end{tabular}}

\caption{Comparisons of objective values of all the methods for solving the sparse regularized least squares problem. The $1^{st}$, $2^{nd}$, and $3^{rd}$ best results are colored with \cone{red}, \ctwo{blue} and \cthree{green}, respectively.}
\label{tab:22}
\vspace{-0.7cm}
\end{table}

\bbb{Data Sets.} Four types of data sets for $\{\bbb{A},\bbb{b}\}$ are considered in our experiments. \textbf{(i)} `random-m-n': We generate the design matrix as $\bbb{A} = \text{randn}(m,n)$, where $\text{randn}(m,n)$ is a function that returns a standard Gaussian random matrix of size $m\times n$. To generate the sparse original signal $\ddot{\bbb{x}} \in \mathbb{R}^n$, we select a support set of size $100$ uniformly at random and set them to arbitrary number sampled from standard Gaussian distribution, the observation vector is generated via $\bbb{b}=\bbb{A}\ddot{\bbb{x}}+\bbb{o}$ with $\bbb{o} = 10\times\text{randn}(m,1)$. \textbf{(ii)} `e2006-m-n': We use the real-world data set `e2006' \footnote{\url{https://www.csie.ntu.edu.tw/~cjlin/libsvmtools/datasets/}} which has been used in sparse optimization \cite{zhang2016accelerated}. We uniformly select $m$ examples and $n$ dimensions from the original data set. \textbf{(iii)} `random-m-n-C': To verify the robustness of \textbf{DEC}, we generate design matrices containing outliers by $\P(\bbb{A})$. Here, $\P(\bbb{A})\in\mathbb{R}^{m\times n}$ is a noisy version of $\bbb{A}\in\mathbb{R}^{m\times n}$ where $2\%$ of the entries of $\bbb{A}$ are corrupted uniformly by scaling the original values by 100 times \footnote{Matlab script: I=randperm(m*n,round(0.02*m*n)); A(I)=A(I)*100.}. We use the same sampling strategy to generate $\bbb{A}$ as in `random-m-n'. Note that the Hessian matrix can be ill-conditioned. \textbf{(iv)} `e2006-m-n-C': We use the same corrupting strategy to generate the corrupted real-world data as in `random-m-n-C'.


\bbb{$\blacktriangleright$ Sparsity Constrained Least Squares Problem.} We compare \textbf{DEC} with 8 state-of-the-art sparse optimization algorithms. \bbb{(i)} Proximal Gradient Method (PGM) \cite{BaoJQS16}, \bbb{(ii)} Accerlated Proximal Gradient Method (APGM), and \bbb{(iii)} Quadratic Penalty Method (QPM) \cite{LuZ13} are gradient-type methods based on iterative hard thresholding. \bbb{(iv)} Subspace Pursuit (SSP) \cite{dai2009subspace}, \bbb{(v)} Regularized Orthogonal Matching Pursuit (ROMP) \cite{needell2010signal}, \bbb{(vi)} Orthogonal Matching Pursuit (OMP) \cite{tropp2007signal}, and \bbb{(vii)} Compressive Sampling Matched Pursuit (CoSaMP)\cite{needell2009cosamp} are greedy algorithms based on iterative support set detection. We use the Matlab implementation in the `sparsify' toolbox\footnote{\url{http://www.personal.soton.ac.uk/tb1m08/sparsify/sparsify.html}}. We also include the comparison with \bbb{(viii)} Convex $\ell_1$ Approximation Method (CVX-$\ell_1$). We use PGM to solve the convex $\ell_1$ regularized problem, with the regulation parameter being swept over $\lambda = \{2^{-10}, 2^{-8},...,2^{10}\}$. The solution that leads to smallest objective after a hard thresholding projection and re-optimization over the support set is selected. Since the optimal solution is expected to be sparse, we initialize the solutions of \{PGM, APGM, QPM, CVX-$\ell_1$, \textbf{DEC}\} to $10^{-7} \times \text{randn}(n,1)$ and project them to feasible solutions. The initial solution of greedy pursuit methods are initialized to zero points implicitly. We vary $s=\{3,~8,~13,~18,~...,50\}$ on different data sets and show the average results of using 5 random initial points.

First, we show the convergence curve and computational efficiency of \textbf{DEC} by comparing with gradient-type methods \{PGM,~APGM,~QPM\}. Several observations can be drawn from Figure \ref{fig:1}. \bbb{(i)} PGM and APGM achieve similar performance and they get stuck into poor local minima. \bbb{(ii)} \textbf{DEC} is more effective than \{PGM,~APGM\}. In addition, we find that as the parameter $k$ becomes larger, more higher accuracy is achieved. \bbb{(iii)} \textbf{DEC-R$0$G$2$} converges quickly but it generally leads to worse solution quality than \textbf{DEC-R$2$G$0$}. Based on this observation, we conclude that a combined random and greedy strategy is preferred. \bbb{(iv)} \textbf{DEC} generally takes less than 30 seconds to converge in \emph{all} our instances with obtaining reasonably good accuracy.

Second, we show the experimental results on sparsity constrained least squares problems with varying the cardinality $s$. Several conclusions can be drawn from Figure \ref{fig2}. \bbb{(i)} The methods \{PGM,~APGM,~QPM\} based on iterative hard thresholding generally lead to bad performance. \bbb{(ii)} OMP and ROMP are not stable and sometimes they achieve bad accuracy. \bbb{(iii)}  \textbf{DEC} presents comparable performance to the greedy methods on \{`random-256-1024', ` random-256-2048'\} but it significantly and consistently outperforms the greedy methods on the other data sets.

\bbb{$\blacktriangleright$ Sparse Regularized Least Squares Problem.} We use Proximal Gradient Method (PGM) and Accelerated Proximal Gradient Method (APGM) \cite{nesterov2013introductory,beck2009fast} to solve the $\ell_0$ norm problem directly, leading to two compared methods \bbb{(i)} PGM-$\ell_0$ and \bbb{(ii)} APGM-$\ell_0$. We apply PGM to solve the convex $\ell_1$ relaxation and nonconvex $\ell_p$ relaxation of the original $\ell_0$ norm problem, resulting in additional two methods \bbb{(iii)} PGM-$\ell_1$ and \bbb{(iv)} PGM-$\ell_p$. We compare \textbf{DEC} with $\{$PGM-$\ell_0$,~APGM-$\ell_0$,~PGM-$\ell_1$,~PGM-$\ell_p$$\}$. We initialize the solutions of all the methods to $10^{-7} \times \text{randn}(n,1)$. For PGM-$\ell_p$, we set $p=\frac{1}{2}$ and use the efficient closed-form solver \cite{xu2012l} to compute the proximal operator.

We show the objective values of all the methods with varying the hyper-parameter $\lambda$ on different data sets in Table \ref{tab:22}. Two observations can be drawn. \bbb{(i)} PGM-$\ell_p$ achieves better performance than the convex method PGM-$\ell_1$. \bbb{(ii)} \textbf{DEC} generally outperforms the other methods in all our data sets.

We demonstrate the average computing time for the compared methods in Table \ref{tab:tab4}. We have two observations. \bbb{(i)} \textbf{DEC} takes several times longer to converge than the compared methods. \bbb{(ii)} \textbf{DEC} generally takes less than 70 seconds to converge in \emph{all} our instances. We argue that the computational time is acceptable and pays off as \textbf{DEC} achieves significantly higher accuracy. The main bottleneck of computation is on solving the small-sized subproblems using sub-exponential time $\mathcal{O}(2^k)$. The parameter $k$ can be viewed as a tuning parameter to balance the efficacy and efficiency. 



\begin{table}[!h]
\scalebox{0.94}{\begin{tabular}{|c|c|c|c|c|c|}
\hline
&{\tiny PGM-$\ell_0$} & {\tiny APGM-$\ell_0$} &  {\tiny PGM-$\ell_1$} & {\tiny PGM-$\ell_p$ }&{\tiny DEC-R10G2} \\
\hline
r.-256-1024 & $12\pm3$ & $13\pm3$& $5\pm3$& $15\pm3$ & $36\pm3$  \\
r.-256-2048 & $11\pm3$ & $11\pm3$&  $9\pm3$& $16\pm3$& $66\pm7$ \\
e.-5000-1024 & $12\pm3$ & $11\pm3$& $8\pm3$& $14\pm3$ & $45\pm3$ \\
e.-5000-2048 & $12\pm3$ & $10\pm3$ & $12\pm3$ & $5\pm3$ & $65\pm8$ \\
 \hline
\end{tabular}}
\caption{Comparisons of average times (in seconds) of all the methods on different data sets for solving the sparse regularized least squares problem.}
\label{tab:tab4}
\vspace{-0.7cm}
\end{table}

\section{Conclusions}

This paper presents an effective and practical method for solving sparse optimization problems. Our approach takes advantage of the effectiveness of the combinatorial search and the efficiency of coordinate descent. We provided rigorous optimality analysis and convergence analysis for the proposed algorithm. Our experiments show that our method achieves state-of-the-art performance. Our block decomposition algorithm has been extended to solve binary optimization problems \cite{yuanhybrid2017} and sparse generalized eigenvalue problems \cite{YuanSZ19}.

\vspace{9pt}
\begin{spacing}{0.6}
{ \fontsize{7}{14}\selectfont
\noi \textbf{Acknowledgments.} This work was supported by NSFC (U1911401), Key-Area Research and Development Program of Guangdong Province (2019B121204008), NSFC (61772570, U1811461), Guangzhou Research Project (201902010037), Pearl River S\&T Nova Program of Guangzhou (201806010056), and Guangdong Natural Science Funds for Distinguished Young Scholar (2018B030306025).
}

\end{spacing}



{
\bibliographystyle{ACM-Reference-Format}

\bibliography{my}
}

\clearpage

{ \fontsize{18}{18}\selectfont \textbf{Appendix} }
\vspace{5pt}

The appendix section is organized as follows. Section \ref{sect:useful} presents a useful lemma. Section \ref{sect:global:conv},~\ref{sect:rate:sparse}, and \ref{sect:the3} present respectively the proof of Proposition 1, Theorem 2, and Theorem 3.

\section{A Useful Lemma} \label{sect:useful}
The following lemma is useful in our proof.

\begin{lemma} \label{eq:recursive}
 Assume a nonnegative sequence $\{u^t\}_{t=0}^{\infty}$ satisfies $(u^{t+1})^2 \leq C (u^{t} - u^{t+1})$ for some nonnegative constant $C$. We have:
\beq
u^{t} \leq \fractt{\max(2C, \sqrt{Cu^0})}{t}
\eeq

\begin{proof}

We denote $C_1 \triangleq \max(2C, \sqrt{Cu^0})$. Solving this quadratic inequality, we have:
\beq
u^{t+1} \leq -  \fractt{C}{2} +  \fractt{C}{2} \sqrt{1+\fractt{ 4 u^t}{C}}
\eeq
\noi We now show that $u^{t+1} \leq  \fractt{C_1}{t+1}$, which can be obtained by mathematical induction. (i) When $t=0$, we have $u^{1}  \leq - \fractt{C}{2} + \fractt{C}{2} \sqrt{1+ \fractt{4 u^0}{C} } \leq -\fractt{C}{2} + \fractt{C}{2} (1+\sqrt{\fractt{4u^0}{C}} ) = \fractt{C}{2}\sqrt{\fractt{4u^0}{C}} = \sqrt{Cu^0 } \leq \fractt{C_1}{t+1}$. (ii) When $t\geq 1$, we assume that $u^{t} \leq  \fractt{C_1}{t}$ holds. We derive the following results: $t\geq 1 \Rightarrow \fractt{t+1}{t} \leq 2$ $~\overset{(a)}{\Rightarrow}~  C \fractt{t+1}{t} \leq C_1$ $~\overset{(b)}{\Rightarrow}~ $ $ C (\fractt{1}{t} - \fractt{1}{t+1} ) \leq \fractt{C_1}{(t+1)^2}$ $\Rightarrow \fractt{C}{t} \leq    \fractt{C}{t+1} + \fractt{C_1}{(t+1)^2}$$\Rightarrow \fractt{C C_1}{t} \leq    \fractt{C C_1 }{t+1} + \fractt{C^2_1}{(t+1)^2}$$\Rightarrow \fractt{C^2}{4} + \fractt{C C_1}{t} \leq    \fractt{C C_1 }{t+1} + \fractt{C^2_1}{(t+1)^2}+\fractt{C^2}{4}$$\Rightarrow \fractt{C^2}{4} (1+ \fractt{4 C_1 }{Ct}  )  \leq   (\fractt{C}{2} + \fractt{C_1}{t+1})^2$$\Rightarrow \fractt{C}{2} \sqrt{1+ \fractt{4 C_1 }{Ct} }  \leq   \fractt{C}{2} + \fractt{C_1}{t+1}$$\Rightarrow - \fractt{C}{2} +  \fractt{C}{2} \sqrt{1+ \fractt{4 C_1 }{Ct} }  \leq    \fractt{C_1}{t+1}$$~\overset{(c)}{\Rightarrow}~ - \fractt{C}{2} +  \fractt{C}{2} \sqrt{1+ \fractt{4 u^t }{C} }  \leq    \fractt{C_1}{t+1}$ $\Rightarrow   u^{t+1} \leq \fractt{C_1}{t+1}$. Here, step $(a)$ uses $2C \leq C_1$; step $(b)$ uses $\fractt{1 }{t(t+1)} =  \fractt{1}{t} - \fractt{1}{t+1} $; step $(c)$ uses $u^{t} \leq  \fractt{C_1}{t}$.

\end{proof}

\end{lemma}

\section{Proof of Proposition 1} \label{sect:global:conv}

\begin{proof}

\bbb{(a)} Due to the optimality of $\bbb{x}^{t+1}$, we have: $F(\bbb{x}^{t+1}) + \tfrac{\theta}{2} \|\bbb{x}^{t+1}-\bbb{x}^t\|_2^2 \leq F(\bbb{u})+ \fractt{\theta}{2}\|\bbb{u}-\bbb{x}^t\|_2^2$ for all $\bbb{u}_{\bar{B}} =(\bbb{x}^t)_{\bar{B}}$. Letting $\bbb{u}=\bbb{x}^t$, we obtain the sufficient decrease condition:
\beq \label{eq:suff:dec}
F(\bbb{x}^{t+1})   \leq F(\bbb{x}^t)- \tfrac{\theta}{2}\|\bbb{x}^{t+1}-\bbb{x}^{t}\|^2
\eeq

\noi Taking the expectation of $B$ for the sufficient descent inequality, we have $\textstyle\E[F(\bbb{x}^{t+1}) ] \leq F(\bbb{x}^t) - \E[\tfrac{\theta}{2}\|\bbb{x}^{t+1}-\bbb{x}^t\| ]$. Summing this inequality over $i=0,1,2, ...,t-1$, we have: $\textstyle \tfrac{\theta}{2} \sum_{i=0}^t \E[\| \bbb{x}^{i+1}-\bbb{x}^i \|_2^2 ] \leq F(\bbb{x}^0)- F(\bbb{x}^t).$

\noi Using the fact that $ F(\bar{\bbb{x}}) \leq F(\bbb{x}^t)$, we obtain:
\beq \label{eq:conv:bound:xx}
\textstyle \min_{i=1,...,t} \E[\tfrac{\theta}{2}\| \bbb{x}^{i+1}-\bbb{x}^i \|_2^2] &\leq& \textstyle \tfrac{\theta}{2t} \sum_{i=0}^t \E[\| \bbb{x}^{i+1}-\bbb{x}^i \|_2^2] \nn\\
&\leq&    \tfrac{F(\bbb{x}^0) - F( \bar{\bbb{x}})}{t}
\eeq
\noi Therefore, we have $\lim_{t\rightarrow \infty}\E[\|\bbb{x}^{t+1} - \bbb{x}^t\|] = 0$.

\bbb{(b)} We assume that the stationary point is not a block-$k$ stationary point. In expectation there exists a block of coordinates $B$ such that $\bbb{x}^t \notin\arg\min_{\bbb{z}}~\P(\bbb{z};\bbb{x}^t,B)$ for some $B$, where $\P(\cdot)$ is defined in Definition \ref{def:block:k}. However, according to the fact that $\bbb{x}^t=\bbb{x}^{t+1}$ and subproblem (\ref{eq:subprob}) in Algorithm \ref{algo:main}, we have $\bbb{x}^{t+1} \in\arg\min_{\bbb{z}}~\P(\bbb{z};\bbb{x}^t,B)$. Hence, we have $\bbb{x}^t_{B}\neq \bbb{x}^{t+1}_{B}$. This contradicts with the fact that $\bbb{x}^t=\bbb{x}^{t+1}$ as $t\rightarrow \infty$. We conclude that $\bbb{x}^{t}$ converges to the block-$k$ stationary point.

\end{proof}

\section{Proof of Theorem 2} 
\label{sect:rate:sparse}


\begin{proof}

\bbb{(a)} Note that Algorithm \ref{algo:main} solves problem (\ref{eq:subprob}) in every iteration. Using Proposition \ref{proposition:hierarchy}, we have that the solution $\bbb{x}^{t+1}_B$ is also a $L$-stationary point. Therefore, we have $|\bbb{x}^{t+1}|_i \geq \sqrt{{2\lambda}/{(\theta + L)}}$ for all $\bbb{x}^{t+1}_i\neq 0$. Taking the initial point of $\bbb{x}$ for consideration, we have that:
$$|\bbb{x}_i^{t+1}|\geq \min(|\bbb{x}_i^0|,\sqrt{2\lambda/(\theta + L)}),~\forall~ i=1,~2,...,n.$$
\noi Therefore, we have: $\|\bbb{x}^{t+1}-\bbb{x}^{t}\|_2 \geq \delta$. Taking the expectation of $B$, we have the following results: $\E[\|(\bbb{x}^{t+1}-\bbb{x}^{t})_B\|_2^2] = \fractt{k}{n}\|\bbb{x}^{t+1}-\bbb{x}^{t}\|_2^2 \geq \fractt{k}{n} \delta^2$. Every time the support set of $\bbb{x}$ is changed, the objective value is decreased at least by $\E[\tfrac{\theta}{2}\|\bbb{x}^{t+1}-\bbb{x}^{t}\|^2]\geq \fractt{k \theta\delta^2}{2n} \triangleq D$. Combining with the result in (\ref{eq:conv:bound:xx}), we obtain: $\fractt{[2F(\bbb{x}^0) - 2F(\bar{\bbb{x}})]}{t\theta}  \geq \fractt{\delta^2 k}{n}$. Therefore, the number of iterations is upper bounded by ${\bar{J}}$.

\bbb{(b)} We notice that when the support set is fixed, the original problem reduces to a convex problem. Since the algorithm solves the following subproblem: $\bbb{x}^{t+1} \Leftarrow \arg \min_{\bbb{z}}~f(\bbb{z})  + \fractt{\theta}{2} \|\bbb{z}-\bbb{x}^{t}\|^2,~s.t.~\bbb{z}_{\bar{B}} = \bbb{x}^t_{\bar{B}}$, we have the following optimality condition for $\bbb{x}^{t+1}$:
\beq \label{eq:optimality:sparse}
(\nabla f(\bbb{x}^{t+1}))_{B} + \theta (\bbb{x}^{t+1}-\bbb{x}^{t})_B = \bbb{0},~(\bbb{x}^{t+1})_{\bar{B}} = (\bbb{x}^t)_{\bar{B}},
\eeq

 We now consider the case when $f(\cdot)$ is generally convex. We derive the following inequalities:
\beq
&&\E[F(\bbb{x}^{t+1})] - F(\bar{\bbb{x}})\nn\\
&\overset{(a)}{\leq}& \E[\la \nabla f(\bbb{x}^{t+1}) ,  \bbb{x}^{t+1} - \bar{\bbb{x}} \ra],    \nn\\
&\overset{(b)}{\leq}& \E[\fractt{n}{k} \la (\nabla f(\bbb{x}^{t+1}))_{B} , (\bbb{x}^{t+1} - \bar{\bbb{x}})_{B} \ra], \nn\\
&\overset{(c)}{=}& \E[\fractt{n}{k} \la - \theta (\bbb{x}^{t+1}-\bbb{x}^{t})_B , (\bbb{x}^{t+1} - \bar{\bbb{x}})_{B} \ra]\nn\\
&\overset{(d)}{\leq} & \E[\fractt{n}{k}   \theta \|(\bbb{x}^{t+1}-\bbb{x}^{t})_B\|_2 \cdot \|(\bbb{x}^{t+1} - \bar{\bbb{x}})_{B}\|_2] \nn\\
&\overset{(e)}{=} & \E[\fractt{n}{k}   \theta \|(\bbb{x}^{t+1}-\bbb{x}^{t})\|_2 \cdot \|(\bbb{x}^{t+1} - \bar{\bbb{x}})_{B}\|_2] \nn\\
&\overset{(f)}{\leq} &  \E[{\fractt{n}{k} 2\theta  {\rho}{\sqrt{k}}} \|\bbb{x}^{t+1}-\bbb{x}^{t}\|_2] = \E[\nu\|\bbb{x}^{t+1}-\bbb{x}^{t}\|_2] ~~~~~ \label{eq:cost:to:go}
\eeq
\noi where step $(a)$ uses the convexity of $F(\cdot)$; step $(b)$ uses the fact that each block $B$ is picked randomly with probability $k/n$; step $(c)$ uses the optimality condition of $\bbb{x}^{t+1}$ in (\ref{eq:optimality:sparse}); step $(d)$ uses the Cauchy-Schwarz inequality; step $(e)$ uses $\|(\bbb{x}^{t+1}-\bbb{x}^{t})\|_2 = \|(\bbb{x}^{t+1}-\bbb{x}^{t})_B\|_2$; step $(f)$ uses $\|(\bbb{x}^{t+1} - \bar{\bbb{x}})_B\| \leq {\sqrt{k}}\|(\bbb{x}^{t+1} - \bar{\bbb{x}})_B\|_{\infty}\leq {\sqrt{k}}\|\bbb{x}^{t+1} - \bar{\bbb{x}}\|_{\infty}\leq {\sqrt{k}}(\|\bbb{x}^{t+1}\|_{\infty} + \|\bar{\bbb{x}}\|_{\infty})\leq {2 \rho}{\sqrt{k}}$.

Using the result in (\ref{eq:cost:to:go}) and the sufficient decent condition in (\ref{eq:suff:dec}), we derive the following results:
\beq
\E[F(\bbb{x}^{t+1}) - F(\bar{\bbb{x}})] \leq  \E[\nu \sqrt{ \fractt{2}{\theta} \left(F(\bbb{x}^{t}) - F(\bbb{x}^{t+1})\right)  }] ~~~~~~~~~
\eeq
\noi Denoting $\Delta^{t} \triangleq \E[F(\bbb{x}^{t}) - F(\bar{\bbb{x}})]$ and $C \triangleq \fractt{2 \nu^2}{\theta}$, we have the following inequality:
\beq
(\Delta^{t+1})^2 \leq C (\Delta^{t} - \Delta^{t+1})\nn
\eeq
\noi Combining with Lemma \ref{eq:recursive}, we have:
\beq
\E[F(\bbb{x}^{t}) - F(\bar{\bbb{x}})]  \leq {\max( \fractt{4 \nu^2}{\theta} , \sqrt{  \fractt{2 \nu^2 \Delta^0}{\theta } })}/{t} \nn
\eeq
\noi Therefore, we obtain the upper bound of the number of iterations to converge to a stationary point $\bar{\bbb{x}}$ satisfying $F(\bbb{x}^t)-F(\bar{\bbb{x}}) \leq \epsilon$ with fixing the support set. Combing the upper bound for the number of changes $\bar{J}$ for the support set in \bbb{(a)}, we naturally establish the actual number of iterations for Algorithm \ref{algo:main}.

(\bbb{c}) We now consider the case when $f(\cdot)$ is $\sigma$-strongly convex. We derive the following results:
\beq \label{eq:stronlgy:convex:J2}
&&\E[F(\bbb{x}^{t+1}) - F(\bar{\bbb{x}})   ] \nn\\
&\overset{(a)}{\leq}& \E[-\fractt{\sigma}{2} \|\bar{\bbb{x}}-\bbb{x}^{t+1}\|_2^2 - \la \bar{\bbb{x}}-\bbb{x}^{t+1}, \nabla f(\bbb{x}^{t+1}) \ra]   \nn\\
&\overset{(b)}{\leq}& \E[\fractt{1}{2\sigma} \| \nabla f(\bbb{x}^{t+1}) \|_2^2] \nn\\
&\overset{(c)}{=}&\E[ \fractt{1}{2\sigma}  \|(\nabla f(\bbb{x}^{t+1}))_{B}  \|_2^2 \times \fractt{n}{k}] \nn\\
&\overset{(d)}{=}& \E[\fractt{1}{2\sigma}  \|   \theta (\bbb{x}^{t+1}-\bbb{x}^{t})_B   \|_2^2 \times  \fractt{n}{k}] \nn\\
&\overset{}{=}& \E[  -\fractt{\theta^2 n}{2\sigma k} \| \bbb{x}^{t+1}-\bbb{x}^{t}\|_2^2]  \nn\\
&\overset{(e)}{\leq}& \E[   \fractt{\theta^2 n}{2 \sigma k} \fractt{2}{\theta} \left(F(\bbb{x}^{t}) - F(\bbb{x}^{t+1})]\right) \nn\\
&\overset{(f)}{=}&\E[  \varpi\left( [F(\bbb{x}^{t}) - F(\bar{\bbb{x}}) ] -  [F(\bbb{x}^{t+1}) - F(\bar{\bbb{x}})]   \right) ]
\eeq
\noi where step $(a)$ uses the strong convexity of $f(\cdot)$; step $(b)$ uses $\forall \bbb{x},\bbb{y},~-\fractt{\sigma}{2}\|\bbb{x}\|_2^2 - \la \bbb{x},\bbb{y}\ra \leq \fractt{1}{2 \sigma}\|\bbb{y}\|_2^2$; step $(c)$ uses $\E[\| \bbb{w}_{B}\|_2^2] = \fractt{k}{n}\|\bbb{w}\|_2^2$; step $(d)$ uses the optimality of $\bbb{x}^{t+1}$; step $(e)$ uses the sufficient condition in $(\ref{eq:suff:dec})$; step $(f)$ uses $\varpi \triangleq\frac{n\theta}{k \sigma}$.

Rearranging terms for (\ref{eq:stronlgy:convex:J2}), we have: $\fractt{\E[F(\bbb{x}^{t+1}) - F(\bar{\bbb{x}})]}{\E[F(\bbb{x}^{t}) - F(\bar{\bbb{x}})]} \leq \fractt{\varpi}{1+\varpi} = \alpha $. Solving the recursive formulation, we obtain:
\beq \label{eq:aaa}
\E[F(\bbb{x}^{t}) - F(\bar{\bbb{x}})] \leq \E[\alpha^t [F(\bbb{x}^{0}) - F(\bar{\bbb{x}})]], \nn
\eeq
\noi and it holds that $t \leq \log_{\alpha} \left(\fractt{F(\bbb{x}^{t}) - F(\bar{\bbb{x}})}{F(\bbb{x}^{0}) - F(\bar{\bbb{x}})}\right)$ in expectation. Using similar techniques as in \bbb{(b)}, we obtain (\ref{eq:convex:boundt2}).

 \end{proof}

\section{Proof of Theorem 3}
\label{sect:the3}

\begin{proof}
\bbb{(a)} We first consider the case when $f(\cdot)$ is generally convex. We denote $Z_t = \{i:\bbb{x}^t_i=0\}$ and known that the $Z_t$ only changes for a finite number of times. We assume that $Z_t$ only changes at $t=c_1,c_2,...,c_{\bar{J}}$ and define $c_0=0$. Therefore, we have:
\beq
Z_{0} = Z_{1}=,...,Z_{-1+c_{1}} \neq Z_{c_{1}}=Z_{1+c_{1}}=Z_{2+c_{1}} = ,\nn\\
...,=Z_{-1+c_{j}} \neq Z_{c_{j}} = ... \neq  Z_{c_{{\bar{J}}}} = ...\nn
\eeq
\noi with $j=1,...,{\bar{J}}$. We denote $\bar{\bbb{x}}^{c_j}$ as the optimal solution of the following optimization problem:
\beq \label{eq:convex:subprob}
\min_{\bbb{x}}~f(\bbb{x}) ,~s.t.~\bbb{x}_{Z_{c_j}} = 0
\eeq
\noi with $1\leq j\leq {\bar{J}}$.

The solution $\bbb{x}^{c_j}$ changes $j$ times, the objective values decrease at least by $j D$, where $D$ is defined in (\ref{eq:DD}). Therefore, we have:
\beq
F(\bbb{x}^{c_j}) \leq F(\bbb{x}^0) - {j}\times D \nn
\eeq
\noi Combing with the fact that $ F(\bar{\bbb{x}}) \leq F(\bar{\bbb{x}}^{c_j})$, we obtain:
\beq \label{eq:J:bound00}
0\leq F(\bbb{x}^{c_{j}}) - F(\bar{\bbb{x}}^{c_j}) \leq F(\bbb{x}^0) - F(\bar{\bbb{x}}) - j \times D
\eeq
\noi We now focus on the intermediate solutions $\bbb{x}_{c_{j-1}},~\bbb{x}_{1+c_{j-1}}$,\\
$...$,~$\bbb{x}_{-1+c_j},\bbb{x}_{c_j}$. Using part \bbb{(b)} in Theorem \ref{theorem:sparse:regularized:1}, we conclude that to obtain an accuracy such that $ F(\bbb{x}^{c_{j}}) - F(\bar{\bbb{x}}^{c_j}) \leq  D$, it takes at most ${ \max\left( \fractt{4 \nu^2}{\theta} , \sqrt{  \fractt{2 \nu^2 (F(\bbb{x}^{c_{j}}) - F(\bar{\bbb{x}}^{c_j}))}{\theta } }\right)   }/{D}$ iterations to converge to $\bar{\bbb{x}}^{c_j}$, that is,
{\beq \label{eq:J:bound1}
 c_j-c_{j-1} &\leq& {\max( \fractt{4 \nu^2}{\theta} , \sqrt{  \fractt{2 \nu^2 (F(\bbb{x}^{c_{j}}) - F(\bar{\bbb{x}}^{c_j}))}{\theta } })}/{D} \nn\\
&\leq& {\max( \fractt{4 \nu^2}{\theta} , \sqrt{  \fractt{2 \nu^2 (F(\bbb{x}^0) - F(\bar{\bbb{x}}) - j \times D)}{\theta } })}/{D}~~~~~~~~~~~~~~~
 \eeq
 }


\noi Summing up the inequality in (\ref{eq:J:bound1}) for $j = 1,2,...,\bar{J}$ and using the fact that $j\geq 1$ and $c_0=0$, we obtain that:
\beq \label{eq:J:total}
c_{{\bar{J}}} \leq    \fractt{\bar{J}}{D}  \times \max( \fractt{4 \nu^2}{\theta} , \sqrt{  \fractt{2 \nu^2 (F(\bbb{x}^0) - F(\bar{\bbb{x}}) - D)}{\theta } }) \nn
\eeq
\noi After $c_{\bar{J}}$ iterations, Algorithm \ref{algo:main} becomes the proximal gradient method applied to the problem as in (\ref{eq:convex:subprob}). Therefore, the total number of iterations for finding a block-$k$ stationary point $\bar{N_1}$ is bounded by:
\beq
\bar{N_1} &\overset{(a)}{\leq}& c_{\bar{J}}  + {\max( \fractt{4 \nu^2}{\theta} , \sqrt{  \fractt{2 \nu^2 [F(\bbb{x}^{c_{\bar{J}}}) - F(\bar{\bbb{x}})]}{\theta } })}/{\epsilon } \nn\\
&\overset{(b)}{\leq}& c_{\bar{J}}  + {\max( \fractt{4 \nu^2}{\theta} , \sqrt{  \fractt{2 \nu^2 [ F(\bbb{x}^0) - F(\bar{\bbb{x}}) - D ]}{\theta } })}/{\epsilon } \nn\\
&\overset{}{=}& ( \fractt{\bar{J}}{D} + \fractt{1}{\epsilon}) \times  \max( \fractt{4 \nu^2}{\theta} , \sqrt{  \fractt{2 \nu^2 (F(\bbb{x}^0) - F(\bar{\bbb{x}}) - D)}{\theta } })   \nn
\eeq
\noi where step $(a)$ uses the fact that the total number of iterations for finding a stationary point after $\bbb{x}_{c_{\bar{J}}}$ is upper bounded by ${\max( \fractt{4 \nu^2}{\theta} , \sqrt{  \fractt{2 \nu^2 [F(\bbb{x}^{c_{\bar{J}}}) - F(\bar{\bbb{x}})]}{\theta } })}/{\epsilon }$; step $(b)$ uses (\ref{eq:J:bound00}) and $j\geq 1$.

\bbb{(b)} We now discuss the case when $f(\cdot)$ is strongly convex. Using part \bbb{(c)} in Theorem \ref{theorem:sparse:regularized:1}, we have:
\beq\label{eq:J:bound2}
c_j-c_{j-1} \leq \log_{\alpha} \fractt{D}{F(\bbb{x}^0)-F(\bar{\bbb{x}})} ,
\eeq
\noi Summing up the inequality (\ref{eq:J:bound2}) for $j = 1,2,...,\bar{J}$, we obtain:
\beq \label{eq:jjjjjjjjjj}
c_{{\bar{J}}} \leq \log_{\alpha} ( \fractt{D^{\bar{J}} }{(F(\bbb{x}^0)-F(\bar{\bbb{x}}))^{\bar{J}}} ) = \bar{J} \log_{\alpha} ( \fractt{D}{(F(\bbb{x}^0)-F(\bar{\bbb{x}}))} ) \nn
\eeq
\noi Therefore, the total number of iterations $\bar{N_2}$ is bounded by:
{\beq
\bar{N_{2}} &\overset{(a)}{\leq}& c_{\bar{J}}  + \log_{\alpha} (\fractt{\epsilon}{F(\bbb{x}^{c_{\bar{J}}}) - F(\bar{\bbb{x}})}) \overset{(b)}{\leq} c_{\bar{J}}  + \log_{\alpha} (\fractt{\epsilon}{F(\bbb{x}^{0}) - D - F(\bar{\bbb{x}})}) \nn
\eeq}
\noi where step $(a)$ uses the fact that the total number of iterations for finding a stationary point after $\bbb{x}_{c_{\bar{J}}}$ is upper bounded by $\log_{\alpha} \left(\fractt{\epsilon}{F(\bbb{x}^{c_{\bar{J}}}) - F(\bar{\bbb{x}})}\right)$; step $(b)$ uses (\ref{eq:J:bound00}) that $0\leq  F(\bbb{x}^0) - F(\bar{\bbb{x}}) -  t \times D$ and $t\geq 1$. 

\end{proof}


\end{document}